\documentclass[11pt]{article}
\usepackage[letterpaper, portrait, margin=1in]{geometry}

\usepackage{graphicx}
\usepackage{xcolor}
\usepackage{float}

\usepackage[round]{natbib}

\usepackage{amsmath}
\usepackage{amssymb}
\usepackage{amsfonts}
\usepackage{amsthm}
\usepackage{mathtools}

\allowdisplaybreaks[1]

\newtheorem{theorem}{Theorem}[section]

\newtheorem{lemma}[theorem]{Lemma}
\newtheorem{remark}{Remark}[section]

\newcommand{\ignore}[1]{}

\def\defi{\vcentcolon=}
\def\diff{\mathop{}\!\mathrm{d}}

\title{{\bf Optimizing Pricing, Repositioning, En-Route Time, and Idle Time in Ride-Hailing Systems}}

\author{
{\bf Anton J. Kleywegt},  {\bf Hongzhang Shao} \\
School of Industrial and Systems Engineering \\
Georgia Institute of Technology \\
Atlanta, Georgia 30332-0205 \\
Email: anton@isye.gatech.edu}

\date{
November 2021
}

\begin{document}

\maketitle

\begin{abstract}
In ride-hailing systems, en-route time refers to the time that elapses from the moment a car is dispatched to pick up a rider until the rider is picked up.
A fundamental phenomenon in ride-hailing systems is that there is a trade-off between en-route time and the time that a car waits for a dispatch.
In short, if cars spend little time idle waiting for a dispatch, then few cars are available when a rider makes a request, and thus the mean distance between a rider and the closest available car is long, which means that en-route time is long.
This phenomenon is of great importance in ride-hailing, because en-route time increases rapidly as the number of idle cars decreases, and every minute that a car spends en-route is one minute less that the car can transport riders.
In spite of this, the existing literature on price optimization for ride-hailing, and on repositioning optimization for ride-hailing, ignores en-route time.
Initial attempts to take this trade-off for the mean en-route time into account when considering price optimization or repositioning optimization all resulted in intractable optimization problems.
Then we found a way to reformulate a simultaneous price and repositioning optimization problem, that takes this trade-off for the {\em distribution} of en-route time into account, as a tractable convex optimization problem.
We show how the optimal solution can be used to construct policies that perform much better in simulations than the policies proposed in previous papers.
\end{abstract}

\newpage


\section{Introduction}\label{sec:introduction}

In ride-hailing systems, available vehicle time is one of the most important resources.
Vehicle time can be partitioned into four types of activities:
\begin{enumerate}
\item
The time used to transport riders from their origins (pickup locations) to their destinations (dropoff locations).
This is the major utility (and revenue) generating activity of ride-hailing systems.
This time is called on-trip time.
\item
The time used to reposition empty from one location to another location.
In most transportation systems, travel demand is not balanced over space over time scales of the order of the duration of a trip.
Imbalance in travel demand can be mitigated somewhat by pricing incentives, but most imbalance in travel demand is accommodated either by repositioning of vehicles without riders or by parking of vehicles.
Under typical costs, repositioning of vehicles is preferred over parking of vehicles, and therefore, although repositioning does not generate revenue, it is an essential activity in the operation of most ride-hailing systems.
\item
The time that elapses from the moment a vehicle is dispatched to pick up a rider until the rider is picked up.
This time is called en-route time.
Typically, the farther a vehicle has to travel from its location when it is dispatched to the rider's pickup location, the longer the en-route time.
\item
The time that a vehicle waits to be dispatched.
This time is called idle time.
\end{enumerate}
Here, en-route time and idle time are both in some sense unproductive uses of vehicles' time, and therefore one would like to minimize en-route time and idle time.

\subsection{Trade-off between En-Route time and Idle Time}

A fundamental phenomenon in ride-hailing systems is that there is a trade-off between en-route time and idle time --- if one of these times is reduced, the other time increases.
In short, if vehicles spend little time waiting idle for a dispatch, then few vehicles are available when a rider makes a request, and thus the mean distance between a rider and the closest available vehicle is long, which means that en-route time is long.
The trade-off between en-route time and idle time was pointed out by \cite{arnott1996taxi}, who showed for a stylized setting with pickup locations and vehicles uniformly distributed over a space without boundaries, that the mean en-route distance is proportional to the inverse square root of the density of idle vehicles.

To facilitate this calculation, consider a large (to ignore boundary effects) homogeneous (a simplifying assumption) urban zone. Assume that (1) idle cars are distributed according to a two-dimensional (spatial) Poisson process with a constant rate, and (2) the closest idle car is dispatched for each ride request. Let $\rho$ denotes the density of idle cars at the time a ride request comes, and let $S$ denotes the distance from this passenger to the closest idle car. Clearly, the probability that $S > s$ for some $s>0$ is the (Poisson) probability there is no idle car in a neighborhood of radius $s$. We have
\begin{align}
\mathbb{P}(S_{i} \leq s)
&= 1 - \mathbb{P}(S_{i} > s)
= 1 - \exp\left(-\omega s^2\rho\right)
\label{eqn:pickup-prob}
\end{align}
where $\omega$ is a constant that depends on the chosen distance metric, for example, for $L_{1}$-distance, $\omega = 2$; for $L_{2}$-distance, $\omega \approx 3.14$; for $L_{\infty}$-distance, $\omega = 4$. The mean distance from a passenger origin to the closest idle car is
\begin{align}
\mathbb{E}[S] = \int_{0}^{\infty} \mathbb{P}(S_{i} > s) \diff s
= \int_{0}^{\infty} \exp\left(-\omega s^2\rho\right) \diff s
= \sqrt{\frac{\pi}{4\omega\rho}}
\label{eqn:pickup-dist}
\end{align}
which decreases as $\rho$ increases.

This phenomenon is of great importance in ride-hailing, because en-route time increases rapidly as the number of idle vehicles decreases, and every minute that a vehicle spends en-route is one minute less that the vehicle can transport riders.
In spite of this, much of the existing literature on price optimization for ride-hailing such as \cite{banerjee2016multi,banerjee2016pricing,bimpikis2019spatial}, and on repositioning optimization for ride-hailing such as \cite{braverman2019empty}, ignores en-route time --- it is implicitly assumed that en-route time is zero whatever number of idle vehicles are available.
Exceptions include \cite{castillo2017surge} and \cite{xu2020supply}, that consider the trade-off in the stylized setting with pickup locations and vehicles uniformly distributed over a space without boundaries mentioned above.


\subsection{A New Model of Ride-Hailing Operations}

We consider a continuous time, infinite horizon Markov decision process model of a ride-hailing system with average cost per unit time objective.
Space is partitioned into zones.
Ride requests are associated with a pickup zone and a dropoff zone.
Ride requests for each origin-destination pair arrive according to a Poisson process.
The state includes information about the number of vehicles on-trip between each origin-destination pair, the number of vehicles repositioning between each origin-destination pair, the number of vehicles en-route in each pickup zone for each en-route time class (explained in more detail below), and the number of vehicles idle in each zone.
Instead of just considering the mean en-route time for each zone (approximately proportional to the inverse square root of the current number of idle vehicles in the zone), we consider the distribution of en-route time for each zone which depends on the current number of idle vehicles in the zone.
The distribution is approximated by a mixture of exponential distributions, with the mixture distribution depending on the number of idle vehicles in the zone.
Not only does this give a more accurate model than models that use only the mean en-route time (or models that assume that en-route time is zero), but as discussed below, it also results in a more tractable fluid optimization model than for models that use only the mean en-route time.

The decisions in our model include both origin-destination pricing decisions, as well as repositioning decisions.
This is in contrast with previous work that considered pricing without repositioning \cite{banerjee2016multi,banerjee2016pricing,castillo2017surge,bimpikis2019spatial}, or repositioning without pricing \cite{braverman2019empty}.

The Markov decision process is intractable, partly due to the large state space.
A widely used approach to develop approximately optimal (and under appropriate conditions, asymptotically optimal) policies, is to formulate and solve a deterministic fluid optimization problem associated with the Markov decision process, and then to use an optimal solution of the deterministic fluid optimization problem to compute a policy for the Markov decision process.
Such a deterministic fluid optimization problem is usually obtained by replacing all random variables by their means, and by allowing discrete variables to take fractional values (hence the name fluid optimization problem).
An issue is that even the deterministic fluid optimization problem obtained in the usual way, is intractable for the Markov decision process described above.
As mentioned above, instead of replacing the en-route time by its mean, we approximated the distribution of the en-route time.
We showed that the resulting stochastic fluid optimization problem can be solved in polynomial time, by solving an associated conic optimization problem.

The solution of the stochastic fluid optimization problem can be used in various ways to compute a policy for the Markov decision process.
A simple alternative is the following static (open loop) policy:
The optimal prices for the stochastic fluid optimization problem are used as static prices in the Markov decision process.
The optimal repositioning flows for the stochastic fluid optimization problem are used to compute repositioning probabilities for the Markov decision process.
Another alternative is the following state-dependent (closed loop) policy:
The optimal prices for the stochastic fluid optimization problem are used as static prices in the Markov decision process.
Periodically, the state of the Markov decision process is used as input to a linear program that determines the optimal repositioning decisions to move the state to the optimal state for the stochastic fluid optimization problem while satisfying flow balance constraints.

We also consider an extension of the models described above that take the following into account.
In practice, vehicles can be (and many are) dispatched while repositioning.
For example, while a vehicle is repositioning from zone~$A$ to zone~$B$, it may pass through zone~$C$ and be dispatched to pickup a rider in zone~$C$.
In that case, it does not complete its planned repositioning move.
As far as we are aware, this is the first work that incorporates this important feature.

The performance of a number of policies are compared.
The static and state-dependent policies described above are based on the stochastic fluid optimization problem $\mathsf{FP}_{1}$ that takes the dependence of the en-route time distribution on the number of idle vehicles into account.
To demonstrate the importance of repositioning, we also present results for a policy based on a similar stochastic fluid optimization problem that takes the dependence of the en-route time distribution on the number of idle vehicles into account, but without repositioning (as suggested in \cite{castillo2017surge}).
To demonstrate the importance of taking the dependence of the en-route time distribution on the number of idle vehicles into account, we present results for a static policy based on a similar fluid optimization problem $\mathsf{FP}_{2}$, that ignores en-route time as in much of the existing literature on pricing and repositioning for ride-hailing systems (the type of policy used in \cite{braverman2019empty}).

Our results show that both policies we developed can push the system towards the optimal steady-state, and give average revenues that are close to optimal. Comparing with the static policy, the system under the state-dependent policy is more stable. The simulation results also showed that, although \ref{eqn:FP2} gives more optimistic objective values than \ref{eqn:FP1}, the static policy derived from \ref{eqn:FP2} does not perform well consistently. The performance of this policy is comparable with ours when the demand is relatively low and balanced, but works badly when demand is high and unbalanced. On the other hand, under a ``pricing only'' policy with no repositioning, the average profit from the system is very low in under all tests.


\subsection{Related Literature}

\underline{\emph{En-Route Time and Idle Car Density:}}  The research on balancing en-route time and idle time dates back as early as \cite{arnott1996taxi}. Recent discussions of this topic can be found in \cite{castillo2017surge} and \cite{xu2020supply}. The idea in these papers is that, with a fixed number of drivers, a higher demand will make the density of idle cars go lower, and thus make the en-route time go higher. As cars spend more time picking up passengers, the shortage in supply will become even worse, resulting in a matching failure in the system. There are other papers that model en-route time as a function if idle car density. For example, Besbes et al. (2021) formulate a queueing model of a ride-sharing system, where the expected en-route time (part of the queue's service time) is given by the number of cars and requesting passengers in a continuously dispersed linear region. Our paper differs from these works in two important way: First, we give a tractable way of modeling the relationship between en-route time and idle time in optimization problems. Second, we model the ride-sharing system as a closed queueing network, while other works above study a single homogeneous area without any network structure.

\underline{\emph{Pricing and Matching Radius:}} Pricing has been a popular approach in controlling ride-sharing systems. Papers such as \cite{banerjee2016multi}, \cite{banerjee2016pricing} and \cite{bimpikis2019spatial} use pricing to influence the willingness of passengers to take rides, and thus to balance the demand and supply among different zones in a ride-sharing network. Our work follows a similar approach, where we set a price for each pair of origin and destination, with an objective of maximizing the overall profit. \cite{castillo2017surge} use pricing as a method to overcome the matching failure. They also suggest a maximum dispatch radius to further lower the demand. \cite{xu2020supply} extend this idea and suggest that the maximum dispatch radius should be adaptive. We extend this idea though pricing. In our work, we set a sequence of dispatch radiuses to separate passengers by en-route times. We then offer different prices to passengers with different en-route times (even with the same origin and destination). This approach allows us to make finer control over the demand and the average en-route time in each zone.

\underline{\emph{Empty Repositioning:}} Although \cite{castillo2017surge} suggest that the ``wild goose chase'' can be solved by surge pricing, our numerical simulation gives a different result. As we have mentioned, a pricing-only policy gives very low average revenues, especially when demand is high and imbalanced. Thus, we allow the repositioning of empty cars in our paper. Like pricing, empty repositioning has also been a popular control lever in literature. Our work is very close to \cite{braverman2019empty}. They also consider a closed queueing network model of ride-sharing systems with multiple regions and a fixed number of cars. They propose a static repositioning policy, such that when a car arrives to a zone, it choose to either stay or reposition to some other zones, according to probabilities set by the platform. The major difference here is that, while our work models the en-route time as a function of idle car density, their work simply assumes there to be no en-route time. Our simulation shows that under the same policy, ignoring en-route time will lead to bad system performance when demand is high and imbalanced. Besides a similar static policy, we also give a state-dependent repositioning policy. With simulation, we show that the system under the state-dependent policy is more stable than under the static policy.

\section{Model}\label{sec:model}

In this section, we introduce the model described in Section \ref{sec:introduction} as a Markov decision process (MDP).


\underline{\emph{Network Primitives:}} We consider a ride-sharing system with $N>0$ cars serving a set $\mathcal{I}$ of zones in a city. Let $\sigma_{i}$ denotes the area of zone $i$ for each $i \in \mathcal{I}$. Travel time from zone $i$ to $j$ follows an exponential distribution with mean $1/\mu_{ij}$ when a passenger is in the car, and follows an exponential distribution with mean $1/\widetilde{\mu}_{ij}$ without a passenger. Consider cars in each zone $i \in \mathcal{I}$. Let $\mathcal{I}_{i}$ be the set of destination zones that these cars can travel to with a passenger, and let $\widetilde{\mathcal{I}}_{i}$ be the set of destination zones that they can travel to without carrying a passenger. In addition, let
\begin{align*}
\mathcal{IJ} &\defi \left\{ (i,j)
		\ : \ i \in \mathcal{I}
		\ , \ j \in \mathcal{I}_{i} \right\} \\
\widetilde{\mathcal{IJ}} &\defi \left\{ (i,j)
		\ : \ i \in \mathcal{I}
		\ , \ j \in \widetilde{\mathcal{I}}_{i} \right\}
\end{align*}
denote the sets of zone pairs that cars can travel in between with and without a passenger. We assume that $\mathcal{IJ} \neq \varnothing$. In addition, we assume that for any two zones $i$ and $j$, there must be a sequence of paths in $\widetilde{\mathcal{IJ}}$ such that empty cars can go from $i$ and eventually get to $j$. In other words, the platform has full control over the allocation of empty cars.

When a passenger requests a ride, the platform identifies the available car that is closest to this passenger, and estimates the time for that car to come. The distribution of en-route time to a passenger depends on the number of available cars in that region. For tractability, we only consider a finite number $K_{i}$ of discrete en-route times for travels starting from each zone $i$. Let $K_{0} = \max_{i} \{K_{i} : i \in \mathcal{I}\}$. Consider distances $0 = \delta_{0}^{(N)} < \delta_{1}^{(N)} < \delta_{2}^{(N)} < \cdots < \delta_{K_{0}}^{(N)}$. For any $1 \leq k \leq K_{0}$, if the distance from the passenger to the closest idle car is between $\delta_{k-1}^{(N)}$ and $\delta_{k}^{(N)}$, then the en-route time is assumed to be an exponential distribution with mean $1/\nu_{k}$. We assume that $0 < 1/\nu_{1} < \cdots < 1/\nu_{K_{0}}$. Consider each $i \in \mathcal{I}$ and $1 \leq k \leq K_{i}$. Let $Q_{ik}^{(N)}(a_{i})$ denotes the probability that the mean en-route time for a passenger is $1 / \nu_{k}$, given there are $a_{i}$ available cars at region $i$. Based on the derivation in Section \ref{sec:introduction}, we have
\begin{align*}
Q_{ik}^{(N)}(a_{i})
= \exp \left(- \frac{\omega (\delta_{k-1}^{(N)})^2 a_{i}}{\sigma_{i}} \right) -
	\exp \left(- \frac{\omega (\delta_{k}^{(N)})^2 a_{i}}{\sigma_{i}} \right)
	\quad, \quad 1 \leq k \leq K_{i}
\end{align*}
On the other hand, if the distance to the closest car is larger than $\delta_{K_{i}}$, then the passenger is told that "no car is available", and leaves the system.

\begin{remark}
Considering congestion, the values of $\delta_{1}^{(N)}, \delta_{2}^{(N)}, \cdots, \delta_{K_{0}}^{(N)}$ may be different for different $N$. In practice, we keep a constant sequence of mean en-route times $1/\nu_{0}, 1/\nu_{1}, \cdots, 1/\nu_{K_{0}}$, and then find the distances $\delta_{1}^{(N)}, \delta_{2}^{(N)}, \cdots, \delta_{K_{0}}^{(N)}$ accordingly.
\end{remark}

Passengers request to travel from zone $i$ to zone $j$ with rate $N\lambda_{ij}>0$. Each passenger then decides whether to take the ride offer or not based on the en-route time and the price. Consider a passenger going from $i$ to $j$, and is offered a ride with mean en-route time $1/\nu_{k}$ and price $x_{ijk}$. Let $P_{ijk}(x_{ijk})$ denotes the probability that the passenger accepts the offer. We let
\begin{align*}
P_{ijk}(x_{ijk}) \ \ \defi \ \ \frac{ \exp(\alpha_{ijk} - \beta_{ij} x_{ijk}) }
		{ 1 + \exp(\alpha_{ijk} - \beta_{ij} x_{ijk}) }
\end{align*}
which is a multinomial logit choice model with parameters $\alpha_{ijk}$ and $\beta_{ij}$. Here we assume that $\alpha_{ij1} > \cdots > \alpha_{ijK_{i}}$ and $\beta_{ij}>0$, which means passengers always prefer lower price and shorter waiting time. If the passenger accepts the offer, then the car is dispatched, and the pickup process starts.


\underline{\emph{State Space and Action Space:}} Consider any time $t \geq 0$. For each $(i,j) \in \mathcal{IJ}$ and $1 \leq k \leq K_{i}$, let $D_{ijk}^{(N)}(t)$ denotes the number of cars that are en-route picking up passengers (who wants to go from $i$ to $j$) with a mean travel time $1/\nu_{k}$. For each $(i,j) \in \mathcal{IJ}$, let $E_{ij}^{(N)}(t)$ denotes the number of empty cars driving from zone $i$ to $j$, and let $F_{ij}^{(N)}(t)$ denotes the number of full cars (cars carrying a passenger) driving from zone $i$ to $j$. For each $i \in \mathcal{I}$, let $A_{i}^{(N)}(t)$ denotes the number of cars that are idle at zone $i$. In addition, let:
\begin{align*}
A^{(N)}(t) \ \ &\defi \ \ \left(A_{i}^{(N)}(t)
		\ : \ i \in \mathcal{I} \right) \\
D^{(N)}(t) \ \ &\defi \ \ \left(D_{ijk}^{(N)}(t)
		\ : \ 1 \leq k \leq K_{i} \ , \ (i,j) \in \mathcal{IJ} \right) \\
E^{(N)}(t) \ \ &\defi \ \ \left(E_{ij}^{(N)}(t)
		\ : \ (i,j) \in \widetilde{\mathcal{IJ}} \right) \\
F^{(N)}(t) \ \ &\defi \ \ \left(F_{ij}^{(N)}(t)
		\ : \ (i,j) \in \mathcal{IJ} \right)
\end{align*}
The state of the system at time $t$ is characterized by $(A^{(N)}(t), D^{(N)}(t), E^{(N)}(t), F^{(N)}(t))$. Let
\begin{align*}
\mathcal{S}^{(N)} \ \ \defi \ \ \Big\{ \
(a, d, e, f) \ : \
& \sum_{i \in \mathcal{I}} a_{i} +
		\sum_{(i,j) \in \mathcal{IJ}} \left(
		f_{ij} + e_{ij} + \sum_{k=1}^{K_{i}} d_{ijk}
		\right) = N \ ; \\
& a_{i} \in \mathbb{Z}_+ \ , \ i \in \mathcal{I} \ , \\
& d_{ijk} \in \mathbb{Z}_+ \ , \ 1 \leq k \leq K_{i} \ , \ (i,j) \in \mathcal{IJ} \ , \\
& e_{ij} \in \mathbb{Z}_+ \ , \ (i,j) \in \widetilde{\mathcal{IJ}} \ , \\
& f_{ij} \in \mathbb{Z}_+ \ , \ (i,j) \in \mathcal{IJ} \ \ \Big\}
\end{align*}
denotes the state space of the MDP. We have
\begin{align*}
(A^{(N)}(t), D^{(N)}(t), E^{(N)}(t), F^{(N)}(t)) \in \mathcal{S}^{(N)}
\end{align*}
For any $t \geq 0$.

Right after each random event, the platform can choose to relocate idle cars between zones. Let $Y_{ij}(t)$ denotes the number of idle cars the platform chooses to relocate from zone $i$ to zone $j$ at the latest event before $t$. Let $Y(t) \defi \left(Y_{ij}(t) \ : \ (i,j) \in \mathcal{IJ} \right)$. We need $Y(t) \in \mathcal{A}_{1}(A(t))$, where
\begin{align*}
\mathcal{A}_{1}(a) \ \ \defi \ \ \left\{ y \in \mathbb{Z}_+^{Z \times Z}
		\ : \ \sum_{j \in \widetilde{\mathcal{I}}_{i}} y_{ij} \leq a_{i}
		\ , \ i \in \mathcal{I} \right\}
\end{align*}
Once the relocation decision is made, the platform updates the prices of different rides. For each $(i,j) \in \mathcal{IJ}$ and $1 \leq k \leq K_{i}$, the platform sets a price $X_{ijk}(t)$ for rides going from $i$ to $j$ with mean en-route time $1/\nu_{k}$. Let $X(t) \defi \left(X_{ijk}(t) \ : \ 1 \leq k \leq K_{i} \ , \ (i,j) \in \mathcal{IJ} \right)$. Let
\begin{align*}
\mathcal{A}_{2} \ \ \defi \ \ \left\{ x \ : \ x_{ijk} \in \mathbb{R}
		\ , \ 1 \leq k \leq K_{i} \ , \ (i,j) \in \mathcal{IJ} \right\}
\end{align*}
We need $X_{ijk}(t) \in \mathcal{A}_{2}$ for any $t \geq 0$. For any state $s = (a,d,e,f) \in \mathcal{S}^{(N)}$, we let
\begin{align*}
\mathcal{A}(s) \ \ = \ \ \mathcal{A}_{1}(a) \times \mathcal{A}_{2}
\end{align*}
denotes the action space of the process.


\underline{\emph{Transition Rates and Probability:}} There are four types of events in this stochastic process:
\begin{enumerate}
		\itemsep0em
		\item A passenger request arrives, and a car is dispatched to pickup the passenger;
		\item A car arrivals at the pickup location, and starts delivering the passenger;
		\item A car drops off a passenger at the destination, and becomes idle;
		\item A car arrivals at a new zone after repositioning empty.
\end{enumerate}
Next we give the transition rates and probabilities of these four types of events at each state $s = (a, d, e, f) \in \mathcal{S}^{(N)}$, with action $(x,y)$.

First, let $\theta_{ijk}^{\text{(dispatch)}(N)}$ denotes the rate that a car is dispatched to pickup a passenger who wants to go from $i$ to $j$, with mean en-route time $1 / \nu_{k}$. We have
\begin{align*}
\theta_{ijk}^{\text{(dispatch)}(N)} (s,x,y)
		& \ \ = \ \ N \lambda_{ij} Q_{ik}^{(N)}
		\left(a_{i} - \sum_{j \in \widetilde{\mathcal{I}}_{i}} y_{ij} \right)
		P_{ijk}\left(x_{ijk}\right) \\
		& \quad (i,j) \in \mathcal{IJ}  \ , \ 1 \leq k \leq K_{i}
\end{align*}
With rate $\theta_{ijk}^{\text{(dispatch)}(N)}$, the system evolves as:
\begin{align*}
\zeta_{ijk}^{\text{(dispatch)}}(s,x,y) \ \ &= \ \ s' \ \ = \ \ (a', d', e', f')
\end{align*}
where
\begin{align*}
a_{i'}' \ \ &= \ \ a_{i'} - \left(\sum_{j'=1}^Z y_{i'j'}\right) - \mathbf{1}[i' = i]
		& \quad i' \neq i \ ; \ i' \in \mathcal{I} \\
d_{i'j'k'}' \ \ &= \ \ d_{i'j'k'} + \mathbf{1}[(i',j',k') = (i,j,k)]
		& \quad 1 \leq k \leq K_{i} \ , \ (i',j') \in \mathcal{IJ} \\
e_{i'j'}' \ \ &= \ \ e_{i'j'} + y_{i'j'}
		& \quad (i',j') \in \widetilde{\mathcal{IJ}} \\
f_{i'j'}' \ \ &= \ \ f_{i'j'}
		& \quad (i',j') \in \mathcal{IJ}
\end{align*}
Second, let $\theta_{ijk}^{\text{(pickup)}}$ denotes the rate that a car in $d_{ijk}$ arrives at the pickup location. We have
\begin{align*}
\theta_{ijk}^{\text{(pickup)}}(s,x,y)
		& \ \ = \ \ \nu_{k} d_{ijk}
		& \quad (i,j) \in \mathcal{IJ}  \ , \ 1 \leq k \leq K_{i}
\end{align*}
With rate $\theta_{ijk}^{\text{(pickup)}}$, the system evolves as
\begin{align*}
\zeta_{ijk}^{\text{(pickup)}}(s,x,y) \ \ &= \ \ s' \ \ = \ \ (a', d', e', f')
\end{align*}
where
\begin{align*}
a_{i'}' \ \ &= \ \ a_{i'} - \left(\sum_{j'=1}^Z y_{i'j'}\right)
		& \quad i' \neq i \ ; \ i' \in \mathcal{I} \\
d_{i'j'k'}' \ \ &= \ \ d_{i'j'k'} - \mathbf{1}[(i',j',k') = (i,j,k)]
		& \quad 1 \leq k \leq K_{i} \ , \ (i',j') \in \mathcal{IJ} \\
e_{i'j'}' \ \ &= \ \ e_{i'j'} + y_{i'j'}
		& \quad (i',j') \in \widetilde{\mathcal{IJ}} \\
f_{i'j'}' \ \ &= \ \ f_{i'j'} + \mathbf{1}[(i',j') = (i,j)]
		& \quad (i',j') \in \mathcal{IJ}
\end{align*}
Third, let $\theta_{ij}^{\text{(dropoff)}}$ denotes the rate that a car in $f_{ij}$ arrives at the drop-off location. We have
\begin{align*}
\theta_{ij}^{\text{(dropoff)}}(s,x,y)
		& \ \ = \ \ \mu_{ij} f_{ij}
		& \quad (i,j) \in \mathcal{IJ}
\end{align*}
With rate $\theta_{ij}^{\text{(dropoff)}}$, the system evolves as
\begin{align*}
\zeta_{ij}^{\text{(dropoff)}}(s,x,y) \ \ &= \ \ s' \ \ = \ \ (a', d', e', f')
\end{align*}
where
\begin{align*}
a_{i'}' \ \ &= \ \ a_{i'} - \left(\sum_{j'=1}^Z y_{i'j'}\right) + \mathbf{1}[i' = i]
		& \quad i' \neq i \ ; \ i' \in \mathcal{I} \\
d_{i'j'k'}' \ \ &= \ \ d_{i'j'k'}
		& \quad 1 \leq k \leq K_{i} \ , \ (i',j') \in \mathcal{IJ} \\
e_{i'j'}' \ \ &= \ \ e_{i'j'} + y_{i'j'}
		& \quad (i',j') \in \widetilde{\mathcal{IJ}} \\
f_{i'j'}' \ \ &= \ \ f_{i'j'} - \mathbf{1}[(i',j') = (i,j)]
		& \quad (i',j') \in \mathcal{IJ}
\end{align*}
Fourth, let $\theta_{ij}^{\text{(reposition)}}$ denotes the rate that a car in $e_{ij}$ arrives at the destination. We have
\begin{align*}
\theta_{ij}^{\text{(reposition)}}(s,x,y)
		& \ \ = \ \ \widetilde{\mu}_{ij} e_{ij}
		& \quad (i,j) \in \mathcal{IJ}
\end{align*}
With rate $\theta_{ij}^{\text{(reposition)}}$, the system evolves as
\begin{align*}
\zeta_{ij}^{\text{(reposition)}}(s,x,y) \ \ &= \ \ s' \ \ = \ \ (a', d', e', f')
\end{align*}
where
\begin{align*}
a_{i'}' \ \ &= \ \ a_{i'} - \left(\sum_{j'=1}^Z y_{i'j'}\right) + \mathbf{1}[i' = i]
		& \quad i' \neq i \ ; \ i' \in \mathcal{I} \\
d_{i'j'k'}' \ \ &= \ \ d_{i'j'k'}
		& \quad 1 \leq k \leq K_{i} \ , \ (i',j') \in \mathcal{IJ} \\
e_{i'j'}' \ \ &= \ \ e_{i'j'} + y_{i'j'} - \mathbf{1}[(i',j') = (i,j)]
		& \quad (i',j') \in \widetilde{\mathcal{IJ}} \\
f_{i'j'}' \ \ &= \ \ f_{i'j'}
		& \quad (i',j') \in \mathcal{IJ}
\end{align*}
Finally, let
\begin{align*}
\theta^{(N)}(s,x,y) \ \ = \ \ \sum_{(i,j) \in \mathcal{IJ}} \left(
		\theta_{ij}^{\text{(dropoff)}} + \theta_{ij}^{\text{(reposition)}} + \sum_{k=1}^{K_{i}}
		\left(\theta_{ij}^{\text{(dispatch)}(N)} + \theta_{ij}^{\text{(pickup)}}\right)
		\right)
\end{align*}
be the overall transition rate of the system at state $s$, with action $(x,y)$.

In addition, for each pair of $s, s' \in \mathcal{S}^{(N)}$, we let $m(s' | s,x,y)$ dentoes the transition probability from $s$ to $s'$. We do not give the explicit form of $m(s' | s,x,y)$ for conciseness, as transition probabilities can be specified by transition rates.


\emph{\underline{Dynamic Program Formulation:}} We consider the dynamic pricing problem over an infinite horizon. For each state $s \in \mathcal{S}^{(N)}$ and action $(x,y) \in \mathcal{A}(s)$, we define the reward function as
\begin{align*}
g(s,x,y)
\defi \frac{\sum_{(i,j)\in\mathcal{IJ}} \sum_{k=1}^{K_{i}}
		(x_{ijk} - \phi_{ijk})\theta_{ijk}^{\text{(dispatch)}(N)} (s,x,y)}
		{\theta^{(N)}(s,x,y)}
		- \sum_{(i,j)\in\mathcal{IJ}} \psi_{ij} y_{ij}
\end{align*}
where $\phi_{ijk}$ is the cost to deliver a passenger from $i$ to $j$ with en-route time $1/\nu_{k}$, while $\psi_{ij}$ is the cost to reposition a car from $i$ to $j$. We assume that $\phi_{ij1} \leq \cdots \leq \phi_{ijK_{i}}$ for every $i \in \mathcal{I}$.

We aim to find an optimal stationary pricing policy that maximize the long-run average profit for the platform. Given a current state $s \in \mathcal{S}^{(N)}$,  a stationary policy $\pi$ specifies the price array $x(s) = (x_{ijk}(s) \ , \ (i,j) \in \mathcal{IJ} \ , \ 1 \leq k \leq K_{i})$ and the reposition plan $y(s) = (y_{ij}(s)\ , \ (i,j) \in \mathcal{IJ})$. For any stationary policy $\pi$, let $\gamma(\pi)$ denotes the (conservative) long-run average cost under policy $\pi$, given by
\begin{align*}
\gamma(\pi) \defi
		\liminf_{L\to\infty} \frac{\mathbb{E} \left[
		\sum_{l=1}^L X_{ijk}^{-}(T_l) R_{ijk}(T_l) \ \middle| \ s(0) = s(0) \right]}
		{\mathbb{E}[T_L]}
\end{align*}
where $s(0)$ is the initial state. Let $\gamma^* \defi \sup_{\pi} \gamma(\pi)$ denote the optimal (conservative) long-run average cost per period, and let $h^*(s)$ denotes the optimal differential cost in state $s \in \mathcal{S}^{(N)}$. They satisfy the optimality equations
\begin{align*}
h^*(s) = \max_{(x,y) \in \mathcal{A}(s)} \left\{
		g(s,x,y) - \frac{\gamma^*}{\theta^{(N)}(s,x,y)}
		+ \sum_{s'' \in \mathcal{S}^{(N)}} m(s'' | s,x,y) h^*(s'')
\right\}
\end{align*}
and can be found by solving the system of equations above.



\section{Fluid-based Optimization Problems}\label{sec:optimization}

In this section, we introduce a fluid-based optimization problem \ref{eqn:FP1} that approximates the steady-state performance of the MDP. Although \ref{eqn:FP1} is non-convex, we show that its optimal solution can be found efficiently by solving a conic program. We also present a fluid-based optimization problem \ref{eqn:FP2} that is similar to the optimization problem in \cite{braverman2019empty}, and has no en-route time. We show that \ref{eqn:FP2} and can be solved efficiently as well. Solutions from both models will be used for designing policies, which will be compared theoretically and numerically.


\emph{\underline{Fluid-based Optimization Problem:}} First, we consider the following problem:
\begin{subequations}
\begin{align}
\max_{\substack{\bar{a}, \bar{d}, \bar{e}, \bar{f} \\ x, p, q}} \quad
		& \sum_{(i,j) \in \mathcal{IJ}}
		\sum_{k=1}^{K_{i}} \nu_{k} \bar{d}_{ijk} (x_{ijk} - \phi_{ijk})
		\nonumber \\
		& - \sum_{(i,j)\in\widetilde{\mathcal{IJ}}} \psi_{ij} \widetilde{\mu}_{ji} e_{ij}
		\tag{$\mathsf{FP}_{1}$}
		\label{eqn:FP1} \\
\text{s.t.} \quad
&	\sum_{k'=1}^k q_{ik'} \ \ = \ \
		1 - \exp \left(- \frac{\omega \delta_{k}^2 \bar{a}_{i}}{\sigma_{i}} \right)
		& \quad 1 \leq k \leq K_{i} \ , \ i \in \mathcal{I} \label{eqn:FP1-pickup} \\
& p_{ijk} \ \ = \ \ \frac{ \exp(\alpha_{ijk} - \beta_{ij} x_{ijk}) }
		{ 1 + \exp(\alpha_{ijk} - \beta_{ij} x_{ijk}) }
		& \quad 1 \leq k \leq K_{i} \ , \ (i,j) \in \mathcal{IJ} \label{eqn:FP1-choice} \\
& \lambda_{ij} q_{ik}
		p_{ijk} \ \ = \ \ \nu_{k} \bar{d}_{ijk}
		& \quad 1 \leq k \leq K_{i} \ , \ (i,j) \in \mathcal{IJ} \label{eqn:FP1-demand} \\
& \sum_{k=1}^{K_{i}} \nu_{k} \bar{d}_{ijk}
		\ \ = \ \ \mu_{ij} \bar{f}_{ij}
		& (i,j) \in \mathcal{IJ} \label{eqn:FP1-deliver} \\
& \sum_{j \in \widetilde{\mathcal{J}}_{i}} \widetilde{\mu}_{ji} \bar{e}_{ji}
		+ \sum_{j \in \mathcal{J}_{i}} \mu_{ji} \bar{f}_{ji}
		\ \ = \ \ \sum_{j \in \widetilde{\mathcal{I}}_{i}} \widetilde{\mu}_{ij} \bar{e}_{ij}
		+ \sum_{j \in \mathcal{I}_{i}} \mu_{ij} \bar{f}_{ij}
		& \quad i \in \mathcal{I} \label{eqn:FP1-balance} \\
& \sum_{i \in \mathcal{I}} \left( \bar{a}_{i}
		+ \sum_{j \in \widetilde{\mathcal{I}}} \bar{e}_{ij}
		+ \sum_{j \in \mathcal{I}} \left(\bar{f}_{ij}
		+ \sum_{k=1}^{K_{i}} \bar{d}_{ijk} \right) \right)
		\ \ = \ \ 1 \\
& \bar{a} \ , \ \bar{d} \ , \ \bar{e} \ , \ \bar{f} \ \ \geq \ \ 0
\end{align}
\end{subequations}
where $\delta_{k} = \delta_{k}^{(N)} / \sqrt{N}$ for each $k \in \{1, \cdots, K_{0}\}$. The inputs of \ref{eqn:FP1} are network primitives $(\lambda_{ij})$, $(\widetilde{\mu}_{ij})$, $(\mu_{ij})$, the pickup parameter $(\nu_{k})$, $(\delta_{k})$, costs $(\phi_{ijk})$, $(\psi_{ij})$, and choice model parameters $(\alpha_{ijk})$, $(\beta_{ij})$. The variables in the optimization problem are $\bar{a}, \bar{d}, \bar{e}, \bar{f}, x, p, q$, where $x = (x_{ijk})$ represents a static pricing policy $X$, $q = (q_{ik})$ represents ratios of car flows in different pickup queues, and $(\bar{a}, \bar{d}, \bar{e}, \bar{f})$ represents the system state. Intuitively, we can think of $(\bar{a}, \bar{d}, \bar{e}, \bar{f})$ as a placeholder for the stationary distribution of network states in the MDP (scaled by $N$). In \ref{eqn:FP1}, we can interpret $\lambda_{ij} q_{ik} p_{ijk}$ as the rate at which rides are initialized from $i$ to $j$ with an expected en-route time $1 / \nu_{k}$. Constraints (\ref{eqn:FP1-pickup}) controls the ratio of car flow $q_{ik}$ falling into this pickup queue. Such a ride has a price of $x_{ijk}$. Constraints (\ref{eqn:FP1-choice}) gives the acceptance rate $p_{ijk}$ of passengers given $x_{ijk}$. Constraints (\ref{eqn:FP1-demand}) and (\ref{eqn:FP1-deliver}) are Little's Laws for pickup queues and deliver queues respectively. Constraint (\ref{eqn:FP1-balance}) say that the total inflow into each zone must equal the total outflow from that zone.


\emph{\underline{Conic Reformulation of \ref{eqn:FP1}:}} The problem above is not a convex optimization problem as it contains a nonlinear constraints (\ref{eqn:FP1-pickup}), (\ref{eqn:FP1-choice}) and (\ref{eqn:FP1-demand}). Next we show that an optimal solution to \ref{eqn:FP1} can be found efficiently. Specifically, it can be found by solving a conic program. Let
\begin{align*}
\mathcal{K}_{\exp} \defi& \ \text{closure} \ \Big\{ (a_{1}, a_{2}, a_{3}) \ : \ a_{3} \leq a_{2} \log (a_{1} / a_{2})
		\ , \ a_{1} > 0  \ , \ a_{2} > 0 \Big\} \\
= & \Big\{ (a_{1}, a_{2}, a_{3}) \ : \ a_{3} \leq a_{2} \log (a_{1} / a_{2}) \ , \ a_{1} > 0  \ , \ a_{2} > 0 \Big\} \cup
		\Big\{ (a_{1}, 0, a_{3}) \ : \ a_{1} \geq 0 \ , \ a_{3} \leq 0 \Big\}
\end{align*}
denotes the exponential cone. Consider the following conic optimization problem:
\begin{subequations}
\begin{align}
\max_{\substack{\bar{a}, \bar{d}, \bar{e}, \bar{f} \\ u, q}} \quad
		& \sum_{(i,j) \in \mathcal{IJ}}
		\sum_{k=1}^{K_{i}} \left( \nu_{k} u_{ijk} - \nu_{k} \phi_{ijk} \bar{d}_{ijk} \right) \nonumber \\
		& - \sum_{(i,j)\in\widetilde{\mathcal{IJ}}} \psi_{ij} \widetilde{\mu}_{ji} e_{ij}
		\tag{$\mathsf{CP}_{1}$}
		\label{eqn:CP1} \\
\text{s.t.} \quad
& \left(1 - 	\sum_{k'=1}^k q_{ik'}, \ 1, \
		- \frac{\omega \delta_{k}^2 \bar{a}_{i}}{\sigma_{i}} \right)
		\ \ \in \ \ \mathcal{K}_{\exp}
		& \quad 1 \leq k \leq K_{i} \ , \ i \in \mathcal{I}
		\label{eqn:CP1-pickup} \\
& \Big( \lambda_{ij} q_{ik} - \nu_{k} \bar{d}_{ijk}, \ \nu_{k} \bar{d}_{ijk}, \nonumber \\
		& \beta_{ij} \nu_{k} u_{ijk} - \alpha_{ijk} \nu_{k} \bar{d}_{ijk} \Big)
		\ \ \in \ \ \mathcal{K}_{\exp}
		& \quad 1 \leq k \leq K_{i} \ , \ (i,j) \in \mathcal{IJ}
		\label{eqn:CP1-choice} \\
& \sum_{k=1}^{K_{i}} \nu_{k} \bar{d}_{ijk}
		\ \ = \ \ \mu_{ij} \bar{f}_{ij}
		& (i,j) \in \mathcal{IJ}
		\label{eqn:CP1-deliver} \\
& \sum_{j \in \widetilde{\mathcal{J}}_{i}} \widetilde{\mu}_{ji} \bar{e}_{ji}
		+ \sum_{j \in \mathcal{J}_{i}} \mu_{ji} \bar{f}_{ji}
		\ \ = \ \ \sum_{j \in \widetilde{\mathcal{I}}_{i}} \widetilde{\mu}_{ij} \bar{e}_{ij}
		+ \sum_{j \in \mathcal{I}_{i}} \mu_{ij} \bar{f}_{ij}
		& \quad i \in \mathcal{I}
		\label{eqn:CP1-balance} \\
& \sum_{i \in \mathcal{I}} \left( \bar{a}_{i}
		+ \sum_{j \in \widetilde{\mathcal{I}}} \bar{e}_{ij}
		+ \sum_{j \in \mathcal{I}} \left(\bar{f}_{ij}
		+ \sum_{k=1}^{K_{i}} \bar{d}_{ijk} \right) \right)
		\ \ = \ \ 1 \\
& \bar{a} \ , \ \bar{d} \ , \ \bar{e} \ , \ \bar{f} \ , \ q \ \ \geq \ \ 0
\end{align}
\end{subequations}
In \ref{eqn:CP1}, conic constraint (\ref{eqn:FP1-pickup}) corresponds to nonlinear constraint (\ref{eqn:FP1-pickup}) in \ref{eqn:FP1}, while conic constraint (\ref{eqn:FP1-choice}) corresponds to nonlinear constraints (\ref{eqn:FP1-choice}) and (\ref{eqn:FP1-demand}) in \ref{eqn:FP1}. We have:

\begin{theorem}
\label{thm:CP1-FP1}
Let $(\bar{a}^*, \bar{d}^*, \bar{e}^*, \bar{f}^*, u^*, q^*)$ denotes an optimal solution to \ref{eqn:CP1}. \noindent For each $1 \leq k \leq K_{i}$ and $(i,j) \in \mathcal{IJ}$, let
\begin{align*}
x_{ijk}^* \ \ &= \ \
		\begin{cases}
 			u_{ijk}^* / \bar{d}_{ijk}^* & \text{ if } \bar{d}_{ijk}^* > 0 \\
 			0 & \text{ otherwise }
	 	\end{cases} \\
p_{ijk}^* \ \ &= \ \
		\frac{ \exp(\alpha_{ijk} - \beta_{ij} x_{ijk}^*) }
		{ 1 + \exp(\alpha_{ijk} - \beta_{ij} x_{ijk}^*) }
\end{align*}
Then $(\bar{a}^*, \bar{d}^*, \bar{e}^*, \bar{f}^*, x^*, p^*, q^*)$ is an optimal solution to \ref{eqn:FP1}.
\end{theorem}

\begin{remark}
	The theorem above is built upon assumptions we made in Section \ref{sec:model}, including:
\begin{align*}
1/\nu_{K_{0}} \ \ > \ \ \cdots
		\ \ > \ \ 1/\nu_{0} \ \ &> \ \ 0 \\
\delta_{K_{0}} \ \ > \ \ \cdots \ \ > \ \ \delta_{1}
		\ \ > \ \ \delta_{0} \ \ &= \ \ 0 \\
\lambda_{ij} \ \ &> \ \ 0
		& \quad (i,j) \in \mathcal{IJ} \\
\beta_{ij} \ \ &> \ \ 0
		& \quad (i,j) \in \mathcal{IJ} \\
\alpha_{ij1} \ \ > \ \ \cdots \ \ &> \ \ \alpha_{ijK_{i}}
		& \quad (i,j) \in \mathcal{IJ} \\
\end{align*}
\end{remark}

\emph{\underline{Optimization Problem without en-route time:}} Next we introduce a fluid-based optimization model \ref{eqn:FP2} that does not model the trade-off between idle time and en-route time:
\begin{subequations}
\begin{align}
\max_{\substack{\bar{a}, \bar{e}, \bar{f} \\ x, p, q}} \quad
		& \sum_{(i,j) \in \mathcal{IJ}} \mu_{ij} \bar{f}_{ij}  (x_{ij} - \phi_{ij0})
		- \sum_{(i,j)\in\widetilde{\mathcal{IJ}}} \psi_{ij} \widetilde{\mu}_{ji} e_{ij}
		\tag{$\mathsf{FP}_{2}$}
		\label{eqn:FP2} \\
\text{s.t.} \quad
&	(1 - q_{i}) \bar{a}_{i} \ \ = \ \ 0
		& \quad i \in \mathcal{I} \label{eqn:FP2-pickup} \\
& p_{ij} \ \ = \ \ \frac{ \exp(\alpha_{ij0} - \beta_{ij} x_{ij}) }
		{ 1 + \exp(\alpha_{ij0} - \beta_{ij} x_{ij}) }
		& \quad 1 \leq k \leq K_{i} \ , \ (i,j) \in \mathcal{IJ} \label{eqn:FP2-choice} \\
& \lambda_{ij} q_{i} p_{ij} \ \ = \ \ \mu_{ij} \bar{f}_{ij}
		& (i,j) \in \mathcal{IJ} \label{eqn:FP2-deliver} \\
& \sum_{j \in \widetilde{\mathcal{J}}_{i}} \widetilde{\mu}_{ji} \bar{e}_{ji}
		+ \sum_{j \in \mathcal{J}_{i}} \mu_{ji} \bar{f}_{ji}
		\ \ = \ \ \sum_{j \in \widetilde{\mathcal{I}}_{i}} \widetilde{\mu}_{ij} \bar{e}_{ij}
		+ \sum_{j \in \mathcal{I}_{i}} \mu_{ij} \bar{f}_{ij}
		& \quad i \in \mathcal{I} \label{eqn:FP2-balance} \\
& \sum_{i \in \mathcal{I}} \left( \bar{a}_{i}
		+ \sum_{j \in \widetilde{\mathcal{I}}} \bar{e}_{ij}
		+ \sum_{j \in \mathcal{I}} \bar{f}_{ij} \right)
		\ \ = \ \ 1 \\
& 0 \ \ \leq \ \ q_{i} \ \ \leq \ \ 1
		& \quad i \in \mathcal{I} \\
& \bar{a} \ , \ \bar{e} \ , \ \bar{f} \ \ \geq \ \ 0
\end{align}
\end{subequations}
\ref{eqn:FP2} is similar to the fluid-based optimization problem discussed in papers such as \cite{braverman2019empty}, except we use a multinomial logit model to capture the price-demand relationship. In \ref{eqn:FP2}, $(\alpha_{ij0})$ and $(\phi_{ij0})$ are parameters such that $\alpha_{ij0} > \alpha_{ij1}$ and $\phi_{ij0} < \phi_{ij1}$ for every $(i,j) \in \mathcal{IJ}$ (i.e. a shorter en-route time is more attractive and less costly). The en-route time constraint (\ref{eqn:FP1-pickup}) is replaced by (\ref{eqn:FP2-pickup}). Intuitively, $q_{i}$ is the ``availability'' at zone $i$ (percentage of requests that can be served), which takes value $1$ as long as $a_{i}>0$, and can be anything between $0$ and $1$ otherwise. The variable $\bar{d}$ disappears, meaning that pickup takes no time.

\ref{eqn:FP2} can also be transformed into a conic program:
\begin{subequations}
\begin{align}
\max_{\substack{\bar{a}, \bar{e}, \bar{f} \\ u, q}} \quad
		& \sum_{(i,j) \in \mathcal{IJ}}
		\left( \mu_{ij} u_{ij} - \mu_{ij} \phi_{ij0} \bar{f}_{ij} \right)
		- \sum_{(i,j)\in\widetilde{\mathcal{IJ}}} \psi_{ij} \widetilde{\mu}_{ji} e_{ij}
		\tag{$\mathsf{CP}_{2}$}
		\label{eqn:CP2} \\
\text{s.t.} \quad
& \Big( \lambda_{ij} q_{i} - \mu_{ij} \bar{f}_{ij}, \ \mu_{ij} \bar{f}_{ij}, \
		\beta_{ij} \mu_{ij} u_{ij} - \alpha_{ij0} \mu_{ij} \bar{f}_{ij} \Big)
		\ \ \in \ \ \mathcal{K}_{\exp}
		& \quad (i,j) \in \mathcal{IJ}
		\label{eqn:CP2-choice} \\
& \sum_{j \in \widetilde{\mathcal{J}}_{i}} \widetilde{\mu}_{ji} \bar{e}_{ji}
		+ \sum_{j \in \mathcal{J}_{i}} \mu_{ji} \bar{f}_{ji}
		\ \ = \ \ \sum_{j \in \widetilde{\mathcal{I}}_{i}} \widetilde{\mu}_{ij} \bar{e}_{ij}
		+ \sum_{j \in \mathcal{I}_{i}} \mu_{ij} \bar{f}_{ij}
		& \quad i \in \mathcal{I}
		\label{eqn:CP2-balance} \\
& \sum_{i \in \mathcal{I}} \left( \bar{a}_{i}
		+ \sum_{j \in \widetilde{\mathcal{I}}} \bar{e}_{ij}
		+ \sum_{j \in \mathcal{I}} \bar{f}_{ij} \right)
		\ \ = \ \ 1 \\
& 0 \ \ \leq \ \ q_{i} \ \ \leq \ \ 1
		& \quad i \in \mathcal{I} \\
& \bar{a} \ , \ \bar{e} \ , \ \bar{f} \ \ \geq \ \ 0
\end{align}
\end{subequations}

We have
\begin{theorem}
\label{thm:CP2-FP2}
(6) Let $(\bar{a}^*, \bar{e}^*, \bar{f}^*, u^*, q^*)$ denotes an optimal solution to \ref{eqn:CP2}. \noindent For each $(i,j) \in \mathcal{IJ}$, let
\begin{align*}
x_{ij}^* \ \ &= \ \
		\begin{cases}
 			u_{ij}^* / \bar{f}_{ij}^* & \text{ if } \bar{f}_{ij}^* > 0 \\
 			0 & \text{ otherwise }
	 	\end{cases} \\
p_{ij}^* \ \ &= \ \
		\frac{ \exp(\alpha_{ij0} - \beta_{ij} x_{ij}^*) }
		{ 1 + \exp(\alpha_{ij0} - \beta_{ij} x_{ij}^*) }
\end{align*}
Then $(\bar{a}^*, \bar{e}^*, \bar{f}^*, x^*, p^*, q^*)$ is an optimal solution to \ref{eqn:FP2}.
\end{theorem}
In Section \ref{sec:numerical}, we will compare the performance of static policies derived from \ref{eqn:FP1} and \ref{eqn:FP2}.



\section{Static Policys}\label{sec:static} 

In Section \ref{sec:model}, we formulated \ref{eqn:FP1} as an approximation to the MDP's optimal steady state. The next question is how to use \ref{eqn:FP1} in controlling the MDP. One simplest approach here is to formulate a static policy based on an optimal solution to \ref{eqn:FP1}. For pricing, this is easy. Once we get a price vector $x^*$ by solving \ref{eqn:FP1}, we can use $x^*$ as a constant pricing action independent from the system state. For repositioning, however, this is not as straightforward.

In fact, for the current MDP, we cannot formulate a meaningful static reposition policy. By the definition of $\mathcal{A}_{1}(a)$, we can only reposition cars when there are idle cars in a zone. As a result, any static reposition policy for the MDP has to be trivial (no repositioning everywhere). Thus, to build a non-trivial static policy, we will need to slightly modify the current MDP.

\emph{\underline{Modified MDP:}} Consider a Markov decision process that is identical to the model introduced in Section \ref{sec:model}, except we only allow cars to reposition upon arriving at a zone. That is, each time a car reaches a zone, it can either choose to reposition to a different zone, or stay idle at the current zone. The repositioning action is randomly chosen according to a table of probabilities. Consider any time $t \geq 0$. For each $i,j \in \widetilde{\mathcal{IJ}}$, we let $\widetilde{Y}_{ij}(t)$ denotes the probability that a car choose to reposition to $j$ after reaching $i$. For each $i \in \mathcal{I}$, we let $\widetilde{Y}_{ii}(t)$ denotes the probability that a car choose to stay idle after reaching $i$. Let $\widetilde{Y}(t) \defi \left(\widetilde{Y}_{ij}(t) \ : \ i \in \mathcal{I}, \ j \in  \{i\} \cup \widetilde{\mathcal{I}}_{i} \right)$, which is updated after each random event. We need $Y(t) \in \widetilde{\mathcal{A}}_{1}$, where
\begin{align*}
\widetilde{\mathcal{A}}_{1} \ \ \defi \ \ \left\{ \widetilde{y} \in \mathbb{R}_+^{Z \times Z}
		\ : \ \widetilde{y}_{ii} + \sum_{j \in \widetilde{\mathcal{I}}_{i}} \widetilde{y}_{ij} = 1
		\ , \ i \in \mathcal{I} \right\}
\end{align*}
The action space of the modified MDP is then
\begin{align*}
\widetilde{\mathcal{A}} \ \ = \ \ \widetilde{\mathcal{A}}_{1} \times \mathcal{A}_{2}
\end{align*}
which is independent from the system state.

Once $\widetilde{Y}(t)$ is chosen, the repositioning action $Y(t)$ will be determined as follow: (1) If the latest event is a dispatching event or pick-up event, then $Y(t) = 0$; (2) If the latest event is a drop-off event or repositioning event at location $i$ (which means $A_{i}(t) \geq 1$), then $Y_{ii}(t) + \sum_{j \in \widetilde{\mathcal{I}}_{i}} Y_{ij}(t) = 1$, and for every $j \in \{i\} \cup \widetilde{\mathcal{I}}_{i}$, the probability of $Y_{ij}(t) = 1$ is given by $\widetilde{Y}_{ij}$; (3) If the latest event is a drop-off event or repositioning event at location $i$, then $Y_{i'\cdot}(t) = 0$ for every $i' \in \mathcal{I}_{i}$ such that $i' \neq i$. The system will then evolve in the same way as described in Section \ref{sec:model}.

In the rest of this section, our discussion will be based on the modified MDP.

\emph{\underline{Static Policy:}} Next we choose static policies for the modified MDP. Let $(\bar{a}^*, \bar{d}^*, \bar{e}^*, \bar{f}^*, x^*, p^*, q^*)$ be an optimal solution of \ref{eqn:FP1}. Let
\begin{align*}
\widetilde{y}_{ij}^* \ \ &= \ \ \frac{\widetilde{\mu}_{ij} \bar{e}_{ij}^*}
		{\sum_{j' \in \widetilde{\mathcal{I}}_{i}} \widetilde{\mu}_{j'i} \bar{e}_{j'i}^*
		+ \sum_{j' \in \mathcal{I}_{i}} \mu_{j'i} \bar{f}_{j'i}^*}
		& \quad i,j \in \widetilde{\mathcal{IJ}} \\
\widetilde{y}_{ii}^* \ \ &= \ \ 1 - \sum_{j \in \widetilde{\mathcal{I}}_{i}} \widetilde{y}_{ij}^*
		& \quad i \in \mathcal{I}
\end{align*}
Then $(x^*, \widetilde{y}^*)$ is the static policy we choose. (Note that, using the same approach, a static repositioning policy can be chosen based on an optimal solution of \ref{eqn:FP2}.)

\begin{remark}
One possible concern people may have about this repositioning policy, is that cars may reposition repeatedly. That is, a car may choose to reposition again right after finishing a previous repositioning. In fact, this is not a problem. For an optimal solution to \ref{eqn:FP1}, the repositioning car flow will not contain any loop. In other words, no car will be doing a sequence of repositioning and then get back to its starting location. (Otherwise, we are ``wasting'' the resources of cars, which means a better solution should exist.) In this case, even a car may do multiple connected repositioning, they will simply be part of a bigger repositioning plan. We will extend this idea in Section \ref{sec:extension}, where we assume that repositioning can only happen between neighboring zones.
\end{remark}


\section{A State-dependent Repositioning Policy}\label{sec:state-policy}

In this section, we introduce a state-dependent policy that is built upon the solution of \ref{eqn:FP1}, and can be calculated efficiently as well. In this policy, the price vector is still static, and is given by the optimal solution of \ref{eqn:FP1}. Meanwhile, the relocation action depends on the current state, and is obtained by solving a linear optimization problem.

\emph{\underline{Approximated Optimization Problem:}} Suppose that the MDP is at a state $s = (a, d, e, f) \in \mathcal{S}^{(N)}$. To find a good action to take, we need to solve the decision problem:
\begin{align*}
h^*(s) = \max_{(x,y) \in \mathcal{A}(s)} \left\{
		g(s,x,y) - \frac{\gamma^*}{\theta^{(N)}(s,x,y)}
		+ \sum_{s'' \in \mathcal{S}^{(N)}} m(s'' | s,x,y) h^*(s'')
\right\}
\end{align*}
Let $(\bar{a}^*, \bar{d}^*, \bar{e}^*, \bar{f}^*, x^*, p^*, q^*)$ be an optimal solution of \ref{eqn:FP1}. We choose $x^*$ to be our pricing action. The problem becomes
\begin{align*}
\max_{y \in \mathcal{A}_2(s)} \left\{
		g(s,x^*,y) - \frac{\gamma^*}{\theta^{(N)}(s,x^*,y)}
		+ \sum_{s'' \in \mathcal{S}^{(N)}} m(s'' | s,x^*,y) h^*(s'')
\right\}
\end{align*}
Since the state space can be large, it is usually very hard to find $h^*(s'')$ for all states $s'' \in \mathcal{S}^{(N)}$. To solve this optimization problem, we need to approximate $\sum_{s'' \in \mathcal{S}^{(N)}} m(s'' | s,x,y) h^*(s'')$ with a simpler term. We do this in two steps: First, we approximate values of $h^*(s'')$ by a constant $-C < 0$ times some distance measure $L$ between state $s''$ and a ``good state'' $s^* = (a^*, d^*, e^*, f^*)$. In particular, we let
\begin{align*}
(a^*, d^*, e^*, f^*) \ \ &= \ \ N (\bar{a}^*, \bar{d}^*, \bar{e}^*, \bar{f}^*)
\end{align*}
Next, let $s' = (a', d', e', f') \in \mathcal{S}^{(N)}$ be the state right after we take action $y$ in state $s$. Note that $s'$ is deterministic. We have
\begin{subequations}
\begin{align}
a_{i}' \ \ &= \ \ a_{i} - \sum_{j \in \widetilde{\mathcal{I}}_{i}} y_{ij}
		& \quad i \in \mathcal{I} \label{eqn:state-next-1} \\
e_{ij}' \ \ &= \ \ e_{ij} + y_{ij}
		& \quad (i,j) \in \widetilde{\mathcal{IJ}} \label{eqn:state-next-2} \\
d' \ \ &= \ \ d \label{eqn:state-next-3} \\
f' \ \ &= \ \ f \label{eqn:state-next-4}
\end{align}
\end{subequations}
For any possible next state $s'' \in \mathcal{S}^{(N)}$ such that $m(s'' | s,x,y) \neq 0$, the two states $s'$ and $s''$ is only different by one random event (the movement of a single car). Thus, the distance between $s'$ and $s''$ should be small (when $N$ is large), and we can approximate $L(s'',s^*)$ with $L(s',s^*)$. With these two steps, the decision problem becomes
\begin{align*}
\max_{y \in \mathcal{A}_2(s)} \left\{
		g(s,x^*,y) - \frac{\gamma^*}{\theta^{(N)}(s,x,y)}
		- C \times L(s',s^*)
\right\}
\end{align*}

Although the optimization problem becomes simpler, is still not tractable. To make it easier to solve, we make three additional approximations. First, we relax the decision problem such that $y$ can take continuous values. Second, we replace $\theta^{(N)}(s,x,y)$ by a constant (e.g. by $\theta^{(N)}(s,x,0)$), such that $\gamma^* / \theta^{(N)}(s,x,y)$ becomes a constant that can be removed from the problem. Third, recall that:
\begin{align*}
g(s,x,y)
= \frac{\sum_{(i,j)\in\mathcal{IJ}} \sum_{k=1}^{K_{i}}
		(x_{ijk} - \phi_{ijk})\theta_{ijk}^{\text{(dispatch)}(N)} (s,x,y)}
		{\theta^{(N)}(s,x,y)}
		- \sum_{(i,j)\in\mathcal{IJ}} \psi_{ij} y_{ij}
\end{align*}
We replace $\theta_{ijk}^{\text{(dispatch)}(N)} (s,x,y)$ by a constant as well, such that the first term in $g(s,x,y)$ does not depend on $y$, and can be removed. With these simplifications, The decision problem then becomes
\begin{align*}
\min_{y} \quad
		& C \times L(s',s^*) + \sum_{(i,j)\in\mathcal{IJ}} \psi_{ij} y_{ij}
		\tag{$\mathsf{AP}$}
		\label{eqn:AP} \\
\text{s.t.} \quad
		& a_{i} - \sum_{j \in \widetilde{\mathcal{I}}_{i}} y_{ij} \ \ \geq \ \  0
				& \quad i \in \mathcal{I} \\
		& y \ \ \geq \ \ 0
\end{align*}
and it left to define a good and tractable distance measure $L$.

\begin{remark}
The transition rates only depend on $y$ for some of the events. Although $\theta^{(N)}(s,x,y)$ depends on $y$, the impact of $y$ on $\theta^{(N)}(s,x,y)$ is usually small. Thus, we choose to treat $\theta^{(N)}(s,x,y)$ as a constant for simplicity. On the other hand, $\theta_{ijk}^{\text{(dispatch)}(N)} (s,x,y)$ depend largely on $y$. As an extreme case, we can reposition every available car from every location, such that the dispatching rate becomes zero. Still, the first term in $g(s,x,y)$ is highly non-linear, and the expected revenue from a single event is relatively small comparing with other terms. Thus, we choose to only keep the second term in $g(s,x,y)$. The optimization problem will take a focus on minimizing the distance $L(s',s^*)$. When there are multiple repositioning plans giving similar $L(s',s^*)$, the second term in $g(s,x,y)$ will help us to find the one that is least costly.
\end{remark}

\emph{\underline{Distance between Network States:}} Next, we define $L$ and make the optimization problem explicit. It takes time for cars to travel between zones. Thus, it may not be a good idea to measure distance between $s'$ and $s^*$ directly (e.g. with $L_2$ norm). In this section, we measure this distance by solving a minimum-cost network flow problem. We let:
\begin{subequations}
\begin{align}
L(s',s^*) = \nonumber \\
\min_{\substack{d^{+}, d^{-}, \\ e^{+}, e^{-}, f^{-}}} \quad
		& \sum_{(i,j) \in \mathcal{IJ}} \left(
				\frac{e_{ij}^{+}}{\tau} + \frac{e_{ij}^{-}}{\widetilde{\mu}_{ij}}
				+ \frac{f_{ij}^{-}}{\mu_{ij}} \right) \nonumber \\
		+& \sum_{(i,j) \in \mathcal{IJ}} \sum_{k=1}^{K} \left(
				\frac{d_{ijk}^{+}}{\lambda_{ij} p_{ijk}^* q_{ik}^*}
				+ \frac{d_{ijk}^{-}}{\nu_{k}} \right)
		\tag{$\mathsf{LP}$}
		\label{eqn:LP} \\
\text{s.t.} \quad
		& a_{i}' + \sum_{j \in \mathcal{I}_{i}} \left( f_{ji}^{-} + e_{ji}^{-}
				- e_{ij}^{+} - \sum_{k=1}^{K} d_{ijk}^{+} \right)
				\ \ = \ \ a_{i}^*
				& \quad i \in \mathcal{I} \\
		& d_{ijk}' - d_{ijk}^{-} + d_{ijk}^{+} \ \ = \ \ d_{ijk}^*
				& \quad 1 \leq k \leq K_{i} \ , \ (i,j) \in \mathcal{IJ} \\
		& e_{ij}' - e_{ij}^{-} + e_{ij}^{+} \ \ = \ \ e_{ij}^*
				& \quad (i,j) \in \widetilde{\mathcal{IJ}} \\
		& f_{ij}' - f_{ij}^{-} + \sum_{k=1}^{K} d_{ijk}^{-} \ \ = \ \  f_{ij}^*
				& \quad (i,j) \in \mathcal{IJ} \\
		& d^{+} \ , \ d^{-} \ , \ e^{+} \ , \ e^{-} \ , \ f^{-} \ \ \geq \ \ 0
\end{align}
\end{subequations}
In \ref{eqn:LP}, each queue (variable) in the network states is treated as a node in network. For example, if $a_{i}' > a_{i}^*$, then the corresponding node is a source node; otherwise, it is a sink node. Since the total number of cars is a constant, the network is balanced. Edges in this network are determined by the system dynamic of the MDP. For each $(i,j) \in \mathcal{IJ}$ and $1 \leq k \leq K_{i}$, the variable $d_{ijk}^{+}$ represents the inflow to node $d_{ijk}'$ from an external source (arrival of passengers). Similarly, variable $d_{ijk}^{-}$ represents the flow that leaves node $d_{ijk}'$ and enters node $f_{ij}'$ (pickup of passengers). For each $(i,j) \in \mathcal{IJ}$, the variable $e_{ij}^{+}$ represents the flow from node $a_{i}'$ to node $e_{ij}'$ (dispatch of reposition cars), the variable $e_{ij}^{-}$ represents the flow from node $e_{ij}'$ to node $a_{j}'$ (arrival of reposition cars), while the variable $f_{ij}^{-}$ represents the flow from node $f_{ij}'$ to node $a_{j}'$ (drop-off of passengers). The objective of \ref{eqn:LP} is to minimize the total car time needed in transforming $(a', d', e', f')$ into $(a^*, d^*, e^*, f^*)$ (pretending that we have full control over car flows and passenger arrivals). For $d_{ijk}^{-}$, $e_{ij}^{-}$ and $f_{ij}^{-}$, the time costs per car are simply constant travel times $1 / \nu_{k}$, $1 / \mu_{ij}$ and $1 / \widetilde{\mu}_{ij}$ respectively. For $d_{ijk}^{+}$, the time cost per car is the inverse of passenger arrival rate, which can change as the number of available cars changes. In \ref{eqn:LP}, we choose to approximate it by the arrival rate at the optimal state, which is $\lambda_{ij} p_{ijk}^* q_{ik}^*$. For $e_{ij}^{+}$, the cost can either be the time until the next event (when we can dispatch cars for repositioning again), or longer (if we do not have enough idle cars). the value of parameter $\tau$ is something we can choose base on this fact.

After substituting $L(s',s^*)$ by \ref{eqn:LP}, the decision problem \ref{eqn:AP} becomes:
\begin{subequations}
\begin{align}
\min_{\substack{y, \ d^{+}, d^{-}, \\ e^{+}, e^{-}, f^{-}}} \quad
		& \sum_{(i,j) \in \mathcal{IJ}} \psi_{ij} y_{ij} \nonumber \\
		+ & C \sum_{(i,j) \in \mathcal{IJ}} \left(
				\frac{e_{ij}^{+}}{\tau} + \frac{e_{ij}^{-}}{\widetilde{\mu}_{ij}}
				+ \frac{f_{ij}^{-}}{\mu_{ij}} \right) \nonumber \\
		+ & C \sum_{(i,j) \in \mathcal{IJ}} \sum_{k=1}^{K} \left(
				\frac{d_{ijk}^{+}}{\lambda_{ij} p_{ijk}^* q_{ik}^*}
				+ \frac{d_{ijk}^{-}}{\nu_{k}} \right) \nonumber \\
\text{s.t.} \quad
		& a_{i} + \sum_{j \in \mathcal{I}_{i}} \left( f_{ji}^{-} + e_{ji}^{-}
				- y_{ij} - e_{ij}^{+} - \sum_{k=1}^{K} d_{ijk}^{+} \right)
				\ \ = \ \ a_{i}^*
				& \quad i \in \mathcal{I} \\
		& d_{ijk} - d_{ijk}^{-} + d_{ijk}^{+} \ \ = \ \ d_{ijk}^*
				& \quad 1 \leq k \leq K_{i} \ , \ (i,j) \in \mathcal{IJ} \\
		& e_{ij} - e_{ij}^{-} + y_{ij} + e_{ij}^{+} \ \ = \ \ e_{ij}^*
				& \quad (i,j) \in \widetilde{\mathcal{IJ}} \\
		& f_{ij} - f_{ij}^{-} + \sum_{k=1}^{K} d_{ijk}^{-} \ \ = \ \  f_{ij}^*
				& \quad (i,j) \in \mathcal{IJ} \\
		& a_{i} - \sum_{j \in \widetilde{\mathcal{I}}_{i}} y_{ij} \ \ \geq \ \  0
				& \quad i \in \mathcal{I}  \label{eqn:AP-const5} \\
		& y \ , \ d^{+} \ , \ d^{-} \ , \ e^{+} \ , \ e^{-} \ , \ f^{-} \ \ \geq \ \ 0
\end{align}
\end{subequations}
where the network states $(a', d', e', f')$ after action $y$ is replaced by functions of $(a, d, e, f)$ and $y$ according to transformation (\ref{eqn:state-next-1}), (\ref{eqn:state-next-2}) and (\ref{eqn:state-next-3}). Note that the variables $y_{ij}$ and $e_{ij}^{+}$ always appear together in the constraints. Intuitively, $y_{ij}$ is the number of cars we choose to relocate from $i$ to $j$ at the current state, while $e_{ij}^{+}$ is the number of cars we choose to relocate in the future. Thus, we need to choose constraints $\tau$ and $C$ such that $c_{ij} < C/\tau$ for all $(i,j) \in \mathcal{IJ}$. The idea here is that we want to reposition cars at the current state whenever it is possible (as there is no benefit to wait). However, we may not have enough idle cars to be repositioned (constraint (\ref{eqn:AP-const5})), and some repositioning decisions have to wait before the next state. If $c_{ij} \geq C/\tau$, the decision problem may choose to postponing reposition decisions even they can be made right now.

\emph{\underline{Transformation into Action:}} As mentioned, we choose a static price vector $x^*$ to be our pricing action at this current state $s$. The price vector $x^*$ comes from the solution of \ref{eqn:FP1}. We then choose a repositioning action $y$ by solving the decision problem \ref{eqn:AP}. Let $(y^*, (d^{+})^*, (d^{-})^*, (e^{+})^*, (e^{-})^*, (f^{-})^* )$ be an optimal solution to \ref{eqn:AP}. Note that $y^*$ can be not integral. Here, we choose the action $y$ such that $y_{ij} = \lfloor y_{ij}^*\rfloor$ for every $(i,j) \in \mathcal{IJ}$.

\begin{remark}
Even \ref{eqn:AP} is a linear program and can be solved efficiently, it may be not possible to solve it at every event in a large system. On the other hand, we do not actually need to solve it every time. 	Recall that the current state (after taking the action) and the next state is only different by the movement of one car. Suppose that we have done all necessary repositioning at the current state ($(e^{+})^* = 0$). Then at the next state, there should not be much additional thing to do. On the other hand, if we do not have enough idle cars to reposition ($(e_{ij}^{+})^* > 0$ for some $(i,j) \in \mathcal{IJ}$), then we will not be much better off at the next state. Thus, even we solve \ref{eqn:AP} and take actions at a lower frequency (depending on $N$), the system performance should be close to that when we solve \ref{eqn:AP} at every event.
\end{remark}


\section{Extension: Dispatch Repositioning Cars} \label{sec:extension}

In this section, we discuss an extension to the current model, such that repositioning cars are also allowed to be dispatched for passengers.

To facilitate further calculation, we assume that cars can only reposition between neighboring zones. That is, for each zone $i \in \mathcal{I}$, the set $\widetilde{\mathcal{I}}_{i}$ only contains the zones that are adjacent to $i$. Note that making this assumption does not restricts the model. If the platform wants to reposition a car between some far-away zones, then it can do a sequence of repositioning simultaneously between neighboring zones, and achieve the same goal (in a shorter amount of time). For example, consider three zones $i_{1}, i_{2}, i_{3} \in \mathcal{I}$. Suppose that $i_{1}$ and $i_{3}$ are both neighbors of $i_{2}$, but are not adjacent to each other. Then instead of reposition a car from $i_{1}$ to $i_{3}$, the platform can choose to reposition one car from $i_{1}$ to $i_{2}$, and another car from $i_{2}$ to $i_{3}$ at the same time. Since $i_{1}$ is further away from $i_{3}$, the two short repositions takes less time than the one long reposition. This is different from delivering a passenger, as we cannot ``teleport'' a passenger to a different car when the two cars are not in the same zone.

In general, we can set a certain proportion of repositioning cars to be available for the origin zone (e.g. assuming they have not yet across the border), and the rest for the destination zone. Still, if repositioning cars are dispatched for requests in the origin zone, then the number of cars arriving at the destination zone will change. Although we may immediately send another car for the same repositioning, there may be no available car at the origin zone. In other words, allowing such dispatching will change the flow dynamic, and affect the tractability of our fluid model. Thus, it may be more practical to only allow repositioning cars to be dispatched for requests in the destination zone.

Assume that repositioning cars are only available for passenger requests from the destination zone. Then the fluid-based optimization problem becomes:
\begin{subequations}
\begin{align}
\max_{\substack{\bar{a}, \bar{d}, \bar{e}, \bar{f} \\ x, p, q}} \quad
		& \sum_{(i,j) \in \mathcal{IJ}}
		\sum_{k=1}^{K_{i}} \nu_{k} \bar{d}_{ijk} (x_{ijk} - \phi_{ijk})
		\nonumber \\
		& - \sum_{(i,j)\in\widetilde{\mathcal{IJ}}} \psi_{ij} \widetilde{\mu}_{ji} e_{ij}
		\tag{$\mathsf{FP}_{3}$}
		\label{eqn:FP3} \\
\text{s.t.} \quad
&	\sum_{k'=1}^k q_{ik'} \ \ = \ \
		1 - \exp \left(- \frac{\omega \delta_{k}^2 \left(\bar{a}_{i} +
		\sum_{j \in \widetilde{\mathcal{J}}_{i}} \zeta_{ji} \bar{e}_{ji} \right)}{\sigma_{i}} \right)
		& \quad 1 \leq k \leq K_{i} \ , \ i \in \mathcal{I} \label{eqn:FP3-pickup} \\
& p_{ijk} \ \ = \ \ \frac{ \exp(\alpha_{ijk} - \beta_{ij} x_{ijk}) }
		{ 1 + \exp(\alpha_{ijk} - \beta_{ij} x_{ijk}) }
		& \quad 1 \leq k \leq K_{i} \ , \ (i,j) \in \mathcal{IJ} \label{eqn:FP3-choice} \\
& \lambda_{ij} q_{ik}
		p_{ijk} \ \ = \ \ \nu_{k} \bar{d}_{ijk}
		& \quad 1 \leq k \leq K_{i} \ , \ (i,j) \in \mathcal{IJ} \label{eqn:FP3-demand} \\
& \sum_{k=1}^{K_{i}} \nu_{k} \bar{d}_{ijk}
		\ \ = \ \ \mu_{ij} \bar{f}_{ij}
		& (i,j) \in \mathcal{IJ} \label{eqn:FP3-deliver} \\
& \sum_{j \in \widetilde{\mathcal{J}}_{i}} \widetilde{\mu}_{ji} \bar{e}_{ji}
		+ \sum_{j \in \mathcal{J}_{i}} \mu_{ji} \bar{f}_{ji}
		\ \ = \ \ \sum_{j \in \widetilde{\mathcal{I}}_{i}} \widetilde{\mu}_{ij} \bar{e}_{ij}
		+ \sum_{j \in \mathcal{I}_{i}} \mu_{ij} \bar{f}_{ij}
		& \quad i \in \mathcal{I} \label{eqn:FP3-balance} \\
& \sum_{i \in \mathcal{I}} \left( \bar{a}_{i}
		+ \sum_{j \in \widetilde{\mathcal{I}}} \bar{e}_{ij}
		+ \sum_{j \in \mathcal{I}} \left(\bar{f}_{ij}
		+ \sum_{k=1}^{K_{i}} \bar{d}_{ijk} \right) \right)
		\ \ = \ \ 1 \\
& \bar{a} \ , \ \bar{d} \ , \ \bar{e} \ , \ \bar{f} \ \ \geq \ \ 0
\end{align}
\end{subequations}
where $\zeta_{ji}$ are constants in $[0,1]$. We can still write down a conic transformation:
{\small
\begin{subequations}
\begin{align}
\max_{\substack{\bar{a}, \bar{d}, \bar{e}, \bar{f} \\ u, q}} \quad
		& \sum_{(i,j) \in \mathcal{IJ}} \sum_{k=1}^{K_{i}}
		\left( \nu_{k} u_{ijk} - \nu_{k} \phi_{ijk} \bar{d}_{ijk} \right)
		\nonumber \\
		& - \sum_{(i,j)\in\widetilde{\mathcal{IJ}}} \psi_{ij} \widetilde{\mu}_{ji} e_{ij}
		\tag{$\mathsf{CP}_{3}$}
		\label{eqn:CP3} \\
\text{s.t.} \quad
& \left(1 - 	\sum_{k'=1}^k q_{ik'}, \ 1, \
		- \frac{\omega \delta_{k}^2 \left(\bar{a}_{i} +
		\sum_{j \in \widetilde{\mathcal{J}}_{i}} \zeta_{ji} \bar{e}_{ji} \right)}{\sigma_{i}} \right)
		\ \ \in \ \ \mathcal{K}_{\exp}
		& \quad 1 \leq k \leq K_{i} \ , \ i \in \mathcal{I}
		\label{eqn:CP3-pickup} \\
& \Big( \lambda_{ij} q_{ik} - \nu_{k} \bar{d}_{ijk}, \ \nu_{k} \bar{d}_{ijk}, \
		\beta_{ij} \nu_{k} u_{ijk} - \alpha_{ijk} \nu_{k} \bar{d}_{ijk} \Big)
		\ \ \in \ \ \mathcal{K}_{\exp}
		& \quad 1 \leq k \leq K_{i} \ , \ (i,j) \in \mathcal{IJ}
		\label{eqn:CP3-choice} \\
& \sum_{k=1}^{K_{i}} \nu_{k} \bar{d}_{ijk}
		\ \ = \ \ \mu_{ij} \bar{f}_{ij}
		& (i,j) \in \mathcal{IJ}
		\label{eqn:CP3-deliver} \\
& \sum_{j \in \widetilde{\mathcal{J}}_{i}} \widetilde{\mu}_{ji} \bar{e}_{ji}
		+ \sum_{j \in \mathcal{J}_{i}} \mu_{ji} \bar{f}_{ji}
		\ \ = \ \ \sum_{j \in \widetilde{\mathcal{I}}_{i}} \widetilde{\mu}_{ij} \bar{e}_{ij}
		+ \sum_{j \in \mathcal{I}_{i}} \mu_{ij} \bar{f}_{ij}
		& \quad i \in \mathcal{I}
		\label{eqn:CP3-balance} \\
& \sum_{i \in \mathcal{I}} \left( \bar{a}_{i}
		+ \sum_{j \in \widetilde{\mathcal{I}}} \bar{e}_{ij}
		+ \sum_{j \in \mathcal{I}} \left(\bar{f}_{ij}
		+ \sum_{k=1}^{K_{i}} \bar{d}_{ijk} \right) \right)
		\ \ = \ \ 1 \\
& \bar{a} \ , \ \bar{d} \ , \ \bar{e} \ , \ \bar{f} \ , \ q \ \ \geq \ \ 0
\end{align}
\end{subequations}
}%
and show that
\begin{theorem}
\label{thm:CP3-FP3}
Let $(\bar{a}^*, \bar{d}^*, \bar{e}^*, \bar{f}^*, u^*, q^*)$ denotes an optimal solution to \ref{eqn:CP3}. \noindent For each $1 \leq k \leq K_{i}$ and $(i,j) \in \mathcal{IJ}$, let
\begin{align*}
x_{ijk}^* \ \ &= \ \
		\begin{cases}
 			u_{ijk}^* / \bar{d}_{ijk}^* & \text{ if } \bar{d}_{ijk}^* > 0 \\
 			0 & \text{ otherwise }
	 	\end{cases} \\
p_{ijk}^* \ \ &= \ \
		\frac{ \exp(\alpha_{ijk} - \beta_{ij} x_{ijk}^*) }
		{ 1 + \exp(\alpha_{ijk} - \beta_{ij} x_{ijk}^*) }
\end{align*}
Then $(\bar{a}^*, \bar{d}^*, \bar{e}^*, \bar{f}^*, x^*, p^*, q^*)$ is an optimal solution to \ref{eqn:FP3}.
\end{theorem}

\begin{remark}
	As we consume repositioning cars for passenger request, the arrival rate of these cars (or equivalently, exiting rate of this queue) becomes larger than $(\mu_{ji}')$. We ignored this change in the model above since it is relatively small, but it can be large in some edge cases.
\end{remark}


\section{Numerical Experiments}\label{sec:numerical}

In this section, we compare the performance of different policies with simulation. We consider a simplified model of a city given by \cite{braverman2019empty}. The city consists of five zones: a downtown area, a midtown area, and three suburban areas. The downtown area represents a central business district, where many people work but few people live. Midtown represents a zones with restaurants and night-life, where people visit after work. The suburb areas are residential zones, and do not have as many entertainment options as midtown. For convenience, we enumerate the three suburban areas, the midtown area and the downtown area as 1, 2, 3, 4, 5, respectively. We consider three instances:
\begin{enumerate}
\item
In the first instance, the city experiences a rush hour as people go home from work. Most of the traffic originates in downtown, and flows into the suburbs as people go home after work. The network parameters for this instance are

{\small
\begin{align*}
(\lambda_{ij}) &= \begin{pmatrix}
0.0648 &  0.0108 &  0. &  0.0324  &  0. \cr
0.0108 &  0.0648 &  0. &  0.0324  &  0. \cr
0.  &  0.  &  0.0756 &  0.0324  &  0. \cr
0.0216 &  0.0216 &  0.0216 &  0.0216  &  0.0216 \cr
0.324  &  0.324  &  0.324 &  0.108 &  0.
\end{pmatrix}
\ , \
(1 / \mu_{ij}) &= \begin{pmatrix}
0.15 & 0.25 & 1.25 & 0.2 & 0.4 \cr
0.25 & 0.10 & 1.1 & 0.1 & 0.3 \cr
1.25 & 1.1 & 0.1 & 1 & 0.65 \cr
0.25 & 0.15 & 1 & 0.15 & 0.25 \cr
0.5 & 0.4 & 0.75 & 0.25 & 0.2]
\end{pmatrix}
\end{align*}
}

\item
In the second instance, most of the traffic is headed into midtown as people go out in the evening to restaurants and for entertainment. The network parameters for this instance are

{\small
\begin{align*}
(\lambda_{ij}) &= \begin{pmatrix}
0.072 & 0. & 0. & 0.648 & 0. \cr
0. & 0.048 & 0. & 0.432 & 0. \cr
0. & 0. & 0.048 & 0.432 & 0. \cr
0.024 & 0.024 & 0.024 & 0.384 & 0.024 \cr
0. & 0. & 0. & 0.108 & 0.012
\end{pmatrix}
\ , \
(1 / \mu_{ij}) &= \begin{pmatrix}
0.15 & 0.25 & 1.25 & 0.2 & 0.4 \cr
0.25 & 0.10 & 1.1 & 0.1 & 0.3 \cr
1.25 & 1.1 & 0.1 & 1 & 0.65 \cr
0.2 & 0.1 & 1 & 0.15 & 0.25 \cr
0.4 & 0.3 & 0.65 & 0.25 & 0.2
\end{pmatrix}
\end{align*}
}

\item
In the third instance, traffic flows mainly from midtown to the suburbs as people go home for the
night. The network parameters for this instance are

{\small
\begin{align*}
(\lambda_{ij}) &= \begin{pmatrix}
0.108 & 0.006 & 0. & 0.006 & 0. \cr
0.006 & 0.108 & 0. & 0.006 & 0. \cr
0. & 0. & 0.108 & 0.012 & 0. \cr
0.396 & 0.396 & 0.396 & 0.066 & 0.066 \cr
0. & 0. & 0. & 0.012 & 0.108
\end{pmatrix}
\ , \
(1 / \mu_{ij}) &= \begin{pmatrix}
0.15 & 0.25 & 1.25 & 0.2 & 0.4 \cr
0.25 & 0.10 & 1.1 & 0.1 & 0.3 \cr
1.25 & 1.1 & 0.1 & 1 & 0.65 \cr
0.2 & 0.1 & 1 & 0.15 & 0.25 \cr
0.4 & 0.3 & 0.65 & 0.25 & 0.2
\end{pmatrix}
\end{align*}
}
\end{enumerate}

For numerical test, we choose $\omega = 4$. For each $k = 1, \cdots, K_0$, we choose $\delta_k = k$ and $1 / \nu_k = k / 12$. (Same as $1 / \mu_{ij}$, travel time are in hours. Thus for en-route time, we only consider an integer time of five minutes.) For each $i = \mathcal{I}$, we choose $\sigma_{i} = 1$. In addition, for each $(i,j) \in \mathcal{IJ}$ and $k = 0, \cdots, K_0$, we choose $\beta_{ij} = 20 / 5$ and $\alpha_{ijk} = (10 - 10/\nu_k) / 5$. (In the MNL model, this means any trip has a base value of 10, a value of 20 for each hour spent in delivery, and a value of -10 for each hour spent in pickup. In addition, there is a random component in the utility of each individual passenger, which follows a Gumbel distribution with scale parameter $5$.) For convenience, we let $(\widetilde{\mu}_{ij}) = (\mu_{ij})$, and let $(\phi_{ijk}) = (\psi_{ij}) = 0$. In other words, we focus on maximizing the revenue.

For each instance, we put 200 cars in the city, and simulate 20,000 events. We take the first 10,000 events as initialization, and collect results only from on the second 10,000 events.


\emph{\underline{Average Revenue:}} For each instance, we solve \ref{eqn:FP1} and find its optimal objective value. Intuitively, this objective value gives an upper bound to the expected revenue of the stochastic process under any policy. We then simulate the stochastic process under three policies: no repositioning, static repositioning (as in Section \ref{sec:state-policy}), and state-dependent repositioning (as in Section \ref{sec:state-policy}). For all three policies, we take a constant pricing action as given by the solution of \ref{eqn:FP1}. (When finding a pricing-only policy, a constraint $\bar{e}=0$ is added to \ref{eqn:FP1} since there are no repositioning.) The key observations here are: (1) Without repositioning, the revenue is low; (2) The static repositioning policy performs nicely, and the simulated revenues are quite close to the optimal objective values of \ref{eqn:FP1}; (3) The state-dependent repositioning policy performs slightly better than the static policy, and the corresponding simulated revenues are even closer to the optimal values. The exact numbers are shown in Table \ref{tab:ori-revenue-1}.

We also solve \ref{eqn:FP2} to find its optimal objective value, and test the static policy that is built upon the solution of \ref{eqn:FP2}. We observe that: (1) The optimal objective values of \ref{eqn:FP2} is larger than those of \ref{eqn:FP1} (i.e. assuming no en-route time leads to more optimistic predictions); (2) The actual performance of the static policy is far worse than what predicted by \ref{eqn:FP2}; (3) With relatively low and balanced demand (e.g. in Instance 2), the performance of the static policy is comparable with policies developed from \ref{eqn:FP1}; (4) With high and inbalanced demand (e.g. in Instance 1 and Instance 3), however, there is a big gap between the static policy and policies from \ref{eqn:FP1} in performance. The exact numbers are also provided in Table \ref{tab:ori-revenue-1}. As a summary, ignoring en-route time (and its trade-off with idle time) leads to sub-optimal solutions.

\begin{table}[H]
\hspace{10pt}
\begin{center}
\begin{tabular}{c | c c | c c }
 \hline
Scenario & Fluid Problem & Optimal Revenue & Repositioning Policy & Simulated Revenue \\
[0.5ex] \hline\hline
1 & \ref{eqn:FP1} (no repositioning) & 3.85 & no repositioning & 3.87 \\
 \hline
 & \ref{eqn:FP1} & 12.88 & static & 12.41 \\
  &  &  & state-dependent & 12.50 \\
 \hline
 & \ref{eqn:FP2} & 15.79 & static & 8.19 \\
 \hline\hline
2 & \ref{eqn:FP1} (no repositioning) & 6.44 & no repositioning & 6.52 \\
\hline
  & \ref{eqn:FP1} & 15.00 & static & 14.34 \\
  &  &  & state-dependent & 14.73 \\
 \hline
 & \ref{eqn:FP2} & 20.72 & static & 14.13 \\
 \hline\hline
3 & \ref{eqn:FP1} (no repositioning) & 4.31 & no repositioning & 4.30 \\
 \hline
  &\ref{eqn:FP1} & 13.59 & static & 13.11 \\
  &  &  & state-dependent & 13.47 \\
 \hline
 & \ref{eqn:FP2} & 16.78 & static & 8.92 \\
 \hline
\end{tabular}
\caption{\label{tab:ori-revenue-1} Simulated average revenues under different policies (right), compared with optimal objective values of \ref{eqn:FP1} and \ref{eqn:FP2} (left) for three instances. There are 200 cars in each instance, and the average revenues are calculated from the second 10,000 events.}
\end{center}
\end{table}

We also do the same test with less number of cars. Table \ref{tab:ori-revenue-2} shows test results where we have 100 cars in each instance. The results are very similar to those in Table \ref{tab:ori-revenue-1}, except the simulated average revenues are slightly smaller. (Note that the number of cars should not be too small, as we are rounding down the number of repositioning cars in the state-dependent policy.)

\begin{table}[H]
\hspace{10pt}
\begin{center}
\begin{tabular}{c | c c | c c }
 \hline
Scenario & Fluid Problem & Optimal Revenue & Repositioning Policy & Simulated Revenue \\
[0.5ex] \hline\hline
1 & \ref{eqn:FP1} (no repositioning) & 3.85 & no repositioning & 3.72 \\
 \hline
  & \ref{eqn:FP1} & 12.88 & static & 12.24 \\
  &  &  & state-dependent & 12.50 \\
 \hline
 & \ref{eqn:FP2} & 15.79 & static & 8.53 \\
 \hline\hline
2 & \ref{eqn:FP1} (no repositioning) & 6.44 & no repositioning & 6.31 \\
 \hline
  & \ref{eqn:FP1} & 15.00 & static & 14.09 \\
  &  &  & state-dependent & 14.28 \\
 \hline
 & \ref{eqn:FP2} & 20.72 & static & 13.91 \\
 \hline\hline
3 & \ref{eqn:FP1} (no repositioning) & 4.31 & no repositioning & 4.12 \\
 \hline
  & \ref{eqn:FP1} & 13.59 & static & 12.53 \\
  &  &  & state-dependent & 12.83 \\
 \hline
 & \ref{eqn:FP2} & 16.78 & static & 8.55 \\
 \hline
\end{tabular}
\caption{\label{tab:ori-revenue-2} Simulated average revenues under different policies (right), compared with optimal objective values of \ref{eqn:FP1} and \ref{eqn:FP2} (left) for three instances. There are 100 cars in each instance.}
\end{center}
\end{table}


\emph{\underline{Stability of State-dependent Policy:}} We can take a closer look at the simulation results from the tests above. Table \ref{tab:ori-avg-1} gives the optimal number of available (idle) cars at each zone from the solution of \ref{eqn:FP1}. It also provides the average number of available cars at each zone from the simulation results under each policy. It is clear that, under both the static and the state-dependent policy, the average number of available cars are close to the optimal numbers given by \ref{eqn:FP1}. We then look at the second-order deviation in number of idle cars (the square-root of mean of squared differences between the simulate and the optimal value.) From the results in Table \ref{tab:ori-std}, we can see that the state-dependent policy gives lower deviation, which means it performs better in stabilizing the system states than the static policy.

\begin{table}[H]
\hspace{10pt}
\begin{center}
\begin{tabular}{c | c c | c c c c c c }
 \hline
Scenario & Fluid Problem &  &  & SOL \\
 &  &  &  & SU1 & SU2 & SU3 & MT & DT \\
\hline\hline
1 & \ref{eqn:FP1} & &  & 6.73 & 6.01 & 8.81 & 7.90 & 24.85 \\
\hline
2 & \ref{eqn:FP1} & &  & 15.29 & 12.08 & 13.90 & 14.00 & 3.33 \\
\hline
3 & \ref{eqn:FP1} & &  & 6.63 & 6.13 & 7.48 & 25.40 & 7.06 \\
 \hline\hline
Scenario & Fluid Problem & Repositioning Policy &  & AVG \\
 &  &  &  & SU1 & SU2 & SU3 & MT & DT \\
\hline\hline
1 & \ref{eqn:FP1} & static &  & 8.35 & 9.80 & 6.25 & 8.49 & 25.45 \\
 &  & state-dependent &  & 6.97 & 6.71 & 8.41 & 7.82 & 22.96 \\
 \hline
2 & \ref{eqn:FP1} & static &  & 18.92 & 17.36 & 12.11 & 16.23 & 4.09 \\
 &  & state-dependent &  & 14.55 & 11.57 & 12.31 & 18.35 & 3.95 \\
 \hline
3 & \ref{eqn:FP1} & static &  & 7.44 & 8.43 & 6.69 & 32.22 & 6.55 \\
 &  & state-dependent &  & 7.11 & 6.41 & 7.97 & 22.79 & 6.39 \\
 \hline
\end{tabular}
\caption{\label{tab:ori-avg-1} The optimal number of idle cars allocated to each zone in the solution of \ref{eqn:FP1} (above), and the average number of idle cars in each zone from simulation under each policy (below). There are 200 cars in each instance, and the average numbers are calculated from the second 10,000 events.}
\end{center}
\end{table}

\begin{table}[H]
\hspace{10pt}
\begin{center}
\begin{tabular}{c | c c | c c c c c c }
 \hline
Scenario & Fluid Problem & Repositioning Policy &  & SSR \\
 &  &  &  & SU1 & SU2 & SU3 & MT & DT \\
[0.5ex] \hline\hline
1 & \ref{eqn:FP1} & static &  & 6.46 & 12.18 & 7.91 & 7.45 & 37.79 \\
 &  & state-dependent &  & 5.09 & 3.92 & 5.23 & 1.75 & 26.10 \\
 \hline\hline
2 & \ref{eqn:FP1} & static &  & 44.34 & 30.01 & 19.55 & 26.29 & 2.83 \\
 &  & state-dependent &  & 6.80 & 3.75 & 6.74 & 27.96 & 2.14 \\
 \hline\hline
3 & \ref{eqn:FP1} & static &  & 4.77 & 6.66 & 4.94 & 96.96 & 5.21 \\
 &  & state-dependent &  & 3.83 & 3.30 & 6.04 & 22.75 & 5.67 \\
 \hline
\end{tabular}
\caption{\label{tab:ori-std} The second-order deviations between simulated number of available (idle) cars and the optimal number of allocated available cars, under two different policies in three instances.}
\end{center}
\end{table}

On the other hand, if we use \ref{eqn:FP2} (which assumes zero en-route time) to form the static policy, then the actual performance of the stochastic process will be very different from that predicted by \ref{eqn:FP2}. See Table \ref{tab:ori-avg-2} for details.

\begin{table}[H]
\hspace{10pt}
\begin{center}
\begin{tabular}{c | c c | c c c c c c }
 \hline
Scenario & Fluid Problem &  &  & SOL \\
 &  &  &  & SU1 & SU2 & SU3 & MT & DT \\
\hline\hline
1 & \ref{eqn:FP2} & &  & 1.64 & 1.64 & 1.64 & 1.64 & 1.64 \\
\hline
2 & \ref{eqn:FP2} & &  & 0.57 & 0.57 & 0.57 & 0.57 & 0.57 \\
\hline
3 & \ref{eqn:FP2} & &  & 1.62 & 1.62 & 1.62 & 1.62 & 1.62 \\
 \hline\hline
Scenario & Fluid Problem & Repositioning Policy &  & AVG \\
 &  &  &  & SU1 & SU2 & SU3 & MT & DT \\
\hline\hline
1 & \ref{eqn:FP2} & static &  & 12.55 & 48.74 & 21.44 & 11.08 & 3.50 \\
 \hline
2 & \ref{eqn:FP2} & static &  & 6.38 & 8.06 & 7.84 & 9.48 & 7.58 \\
 \hline
3 & \ref{eqn:FP2} & static &  & 9.64 & 34.68 & 34.95 & 3.58 & 11.65 \\
 \hline
\end{tabular}
\caption{\label{tab:ori-avg-2} The optimal number of idle cars allocated to each zone in the solution of \ref{eqn:FP2} (above), and the average number of idle cars in each zone from simulation under a static policy built base on \ref{eqn:FP2} (below). }
\end{center}
\end{table}

\begin{figure}[H]
	\centering
		\includegraphics[width=1\columnwidth]{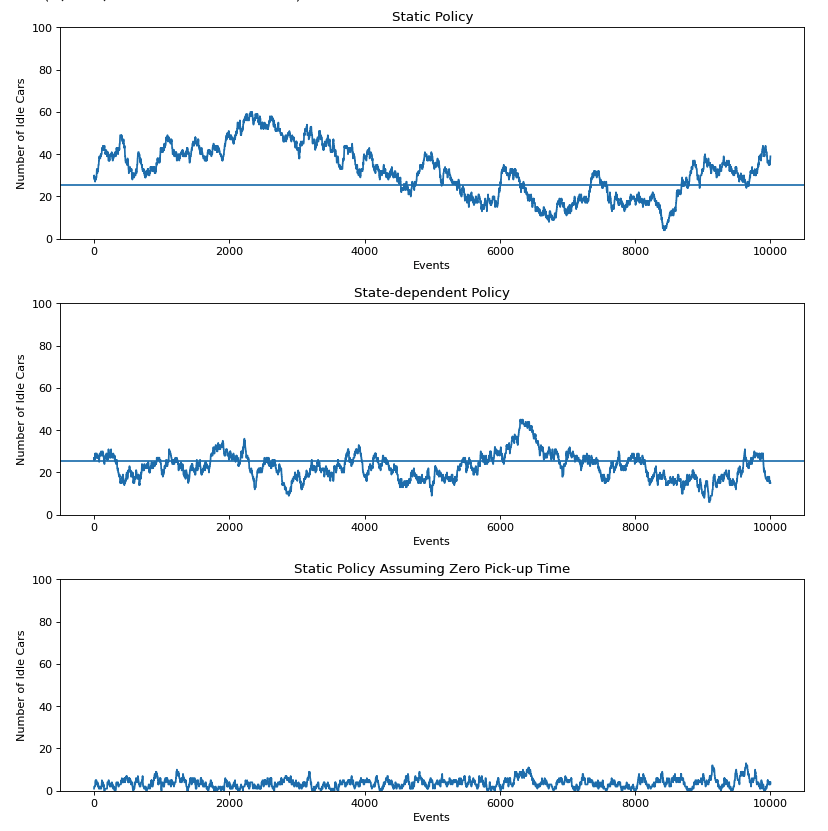} 
		\caption{\label{tab:ori-idle-1} Number of available (idle) cars in the midtown area (the busiest zone) in instance 3, under different policies.}
	\label{plot}
\end{figure}

\begin{figure}[H]
	\centering
		\includegraphics[width=1\columnwidth]{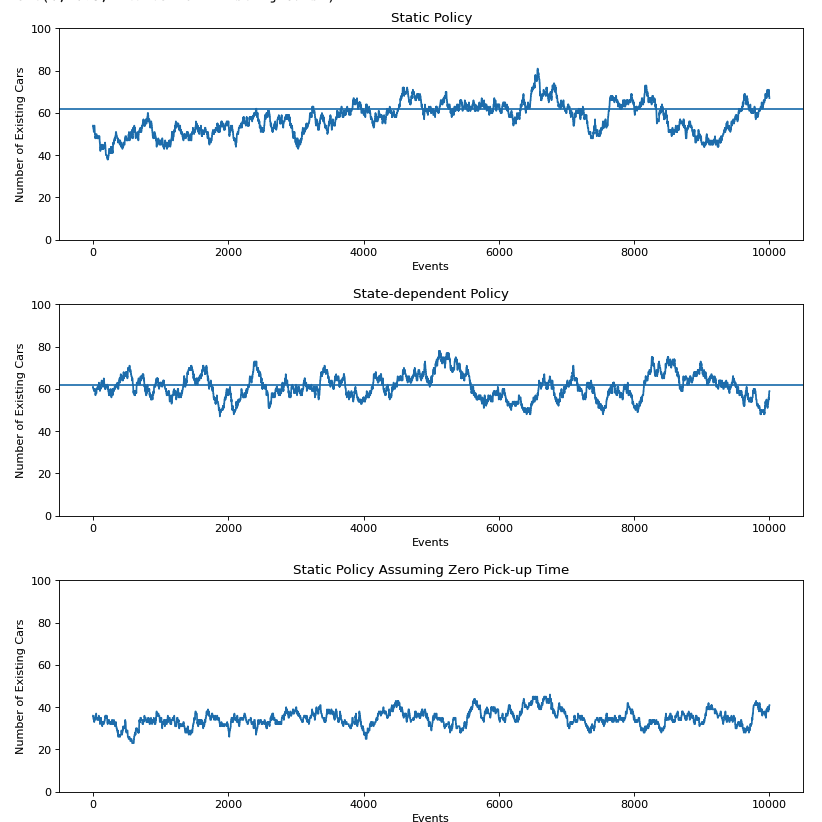} 
		\caption{\label{tab:ori-existing} Number of full cars (cars with a passenger) exiting from the midtown area (the busiest zone) in instance 3, under different policies.}
	\label{plot}
\end{figure}

To make it more visual, Figure \ref{tab:ori-idle-1} shows the number of available cars in the midtown area in instance 3 (during the second 10,000 simulated events) under different policies. Figure \ref{tab:ori-existing} then shows the total number of cars that are flowing out from the midtown area. In this instance, most trips originates from the midtown area, which makes it the busiest zone. As we can see, the number of available cars under the static and state-dependent policy from \ref{eqn:FP1} are both close to that predicted by \ref{eqn:FP1} (around $25$), but the state-dependent policy has a more stable performance. Similar arguments hold for the number of exiting cars. On the other hand, when we use \ref{eqn:FP2} to build the static policy, both the number of available cars and the delivery flow are much smaller, leading to a sub-optimal performance.

\emph{\underline{Expanded Network:}} We repeat the tests above on an ``expanded'' stochastic process and the fluid model \ref{eqn:FP3} as introduced in Section \ref{sec:extension}, where repositioning cars are also available to be dispatched for passenger requests. Table \ref{tab:exp-revenue} and \ref{tab:exp-std} present some of the simulation results. Observations from the previous tests still hold here. For the ``expanded'' stochastic process, the gap in performance between the static and state-dependent policies is even larger, as the state-dependent policy has stronger controls over the repositioning of cars.

\begin{table}[H]
\hspace{10pt}
\begin{center}
\begin{tabular}{c | c c | c c }
 \hline
Scenario & Fluid Problem & Optimal Revenue & Repositioning Policy & Simulated Revenue \\
[0.5ex] \hline\hline
1 & \ref{eqn:FP3} (no repositioning) & 3.85 & no repositioning & 3.81 \\
\hline
  & \ref{eqn:FP3} & 13.97 & static & 12.29 \\
  &  &  & state-dependent & 13.22 \\
 \hline\hline
2 & \ref{eqn:FP3} (no repositioning) & 6.44 & no repositioning & 6.37 \\
\hline
  & \ref{eqn:FP3} & 17.01 & static & 15.11 \\
  &  &  & state-dependent & 15.61 \\
 \hline\hline
3 & \ref{eqn:FP3} (no repositioning) & 4.31 & no repositioning & 4.31 \\
 \hline
  & \ref{eqn:FP3} & 14.66 & static & 10.41 \\
  &  &  & state-dependent & 13.98 \\
 \hline
\end{tabular}
\caption{\label{tab:exp-revenue} Simulated average revenues under different policies (right), compared with optimal objective values of \ref{eqn:FP3} (left) for three instances. There are 200 cars in each instance.}
\end{center}
\end{table}

\begin{table}[H]
\hspace{10pt}
\begin{center}
\begin{tabular}{c | c c | c c c c c c }
 \hline
Scenario & Fluid Problem & Repositioning Policy &  & SSR \\
 &  &  &  & SU1 & SU2 & SU3 & MT & DT \\
[0.5ex] \hline\hline
1 & \ref{eqn:FP3} & static &  & 46.69 & 11.73 & 114.47 & 26.17 & 121.37 \\
 &  & state-dependent &  & 32.5 & 12.68 & 42.18 & 4.72 & 42.84 \\
 \hline\hline
2 & \ref{eqn:FP3} & static &  & 46.03 & 73.90 & 13.87 & 63.14 & 113.52 \\
 &  & state-dependent &  & 26.02 & 13.53 & 13.46 & 30.67 & 101.73 \\
 \hline\hline
3 & \ref{eqn:FP3} & static &  & 24.67 & 1538.12 & 66.79 & 165.21 & 40.29 \\
 &  & state-dependent &  & 8.81 & 3.96 & 22.07 & 38.34 & 8.61 \\
 \hline
\end{tabular}
\caption{\label{tab:exp-std} The second-order deviations between simulated number of available (idle) cars and the optimal number of allocated available cars, under two different policies in three instances. Note that in instance 3, the static policy is not strong enough to ``drag'' the system states back to optimal, and gives a very high deviation in the allocation of idle cars. }
\end{center}
\end{table}

\newpage


\bibliography{library}
\bibliographystyle{plainnat}

\appendix

\newpage


\section{Proofs}


\subsection{Proof of Theorem \ref{thm:CP1-FP1}}
\label{sec:proof-conic}

{\color{purple}
To build the intuition, first note that:
}

{\color{blue}
\begin{lemma}
\label{lem:FP1-CP1}
Consider any $(\bar{a}, \bar{d}, \bar{e}, \bar{f}, x, p, q)$ feasible to \ref{eqn:FP1}. Let
\begin{align*}
u_{ijk} \ \ &= \ \ \bar{d}_{ijk} x_{ijk} \ \
		& \quad 1 \leq k \leq K_{i} \ , \ (i,j) \in \mathcal{IJ}
\end{align*}
Then $(\bar{a}, \bar{d}, \bar{e}, \bar{f}, u, q)$ is feasible to \ref{eqn:CP1}, and has the same objective value as $(\bar{a}, \bar{d}, \bar{e}, \bar{f}, x, p, q)$ in \ref{eqn:FP1}.
\end{lemma}
}%

\begin{proof}
It is clear that the objective value of \ref{eqn:CP1} at $(\bar{a}, \bar{d}, \bar{e}, \bar{f}, u, q)$ is identical to the objective value of \ref{eqn:FP1} at $(\bar{a}, \bar{d}, \bar{e}, \bar{f}, x, p, q)$. We only need to prove that $(\bar{a}, \bar{d}, \bar{e}, \bar{f}, u, q)$ is feasible to \ref{eqn:CP1}. Indeed. Since $(\bar{a}, \bar{d}, \bar{e}, \bar{f}, x, p, q)$ is feasible to \ref{eqn:FP1}, we have
\begin{align*}
&	\sum_{k'=1}^k q_{ik'} \ \ = \ \
		1 - \exp \left(- \frac{\omega \delta_{k}^2 \bar{a}_{i}}{\sigma_{i}} \right) \\
\implies \quad
& - \frac{\omega \delta_{k}^2 \bar{a}_{i}}{\sigma_{i}}
		\ \ = \ \ \ln \left( 1 - \sum_{k'=1}^k q_{ik'} \right) \\
\implies \quad
& - \frac{\omega \delta_{k}^2 \bar{a}_{i}}{\sigma_{i}}
		\ \ \leq \ \ \ln \left( 1 - \sum_{k'=1}^k q_{ik'} \right) \\
\implies \quad
& \left(1 - 	\sum_{k'=1}^k q_{ik'}, \ 1, \
		- \frac{\omega \delta_{k}^2 \bar{a}_{i}}{\sigma_{i}} \right)
		\ \ \in \ \ \mathcal{K}_{\exp}
\end{align*}
For each $1 \leq k \leq K_{i}$ and $i \in \mathcal{I}$. Thus, constraint (\ref{eqn:CP1-pickup}) holds in \ref{eqn:CP1} at $(\bar{a}, \bar{d}, \bar{e}, \bar{f}, u, q)$. Similarly, consider each $1 \leq k \leq K_{i}$ and $(i,j) \in \mathcal{IJ}$. If $\lambda_{ij} q_{ik}>0$, then we have
\begin{align*}
& p_{ijk} \ \ = \ \ \frac{ \exp(\alpha_{ijk} - \beta_{ij} x_{ijk}) }
		{ 1 + \exp(\alpha_{ijk} - \beta_{ij} x_{ijk}) } \\
& \lambda_{ij} q_{ik}
		p_{ijk} \ \ = \ \ \nu_{k} \bar{d}_{ijk} \\
\implies \quad
& \frac{\nu_{k} \bar{d}_{ijk}}{\lambda_{ij} q_{ik}}
		\ \ = \ \ \frac{ \exp(\alpha_{ijk} - \beta_{ij} x_{ijk}) }
		{ 1 + \exp(\alpha_{ijk} - \beta_{ij} x_{ijk}) } \\
\iff \quad
& \frac{\nu_{k} \bar{d}_{ijk}}{\lambda_{ij} q_{ik} - \nu_{k} \bar{d}_{ijk}}
		\ \ = \ \ \exp(\alpha_{ijk} - \beta_{ij} x_{ijk}) \\
\implies \quad
& (\lambda_{ij} q_{ik} - \nu_{k} \bar{d}_{ijk}, \
		\nu_{k} \bar{d}_{ijk}, \ \beta_{ij} \nu_{k} u_{ijk} - \alpha_{ijk} \nu_{k} \bar{d}_{ijk})
		\ \ \in \ \ \mathcal{K}_{\exp}
\end{align*}
On the other hand, if $\lambda_{ij} q_{ik} = 0$, then $\bar{d}_{ijk} = 0$, which implies $u_{ijk} = 0$. We still have
\begin{align*}
& (\lambda_{ij} q_{ik} - \nu_{k} \bar{d}_{ijk}, \
		\nu_{k} \bar{d}_{ijk}, \ \beta_{ij} \nu_{k} u_{ijk} - \alpha_{ijk} \nu_{k} \bar{d}_{ijk})
		\ \ = \ \ (0,0,0)
		\ \ \in \ \ \mathcal{K}_{\exp}
\end{align*}
Thus, constraint (\ref{eqn:CP1-choice}) holds in \ref{eqn:CP1} at $(\bar{a}, \bar{d}, \bar{e}, \bar{f}, u, q)$. Since other constraints in \ref{eqn:CP1} are identical to those in \ref{eqn:FP1}, we can conclude that $(\bar{a}, \bar{d}, \bar{e}, \bar{f}, u, q)$ is feasible to \ref{eqn:CP1}.
\end{proof}

{\color{purple}
On the other hand, we have:
}

{\color{blue}
\begin{lemma}
\label{lem:CP1-FP1}
Consider any $(\bar{a}, \bar{d}, \bar{e}, \bar{f}, u, q)$ feasible to \ref{eqn:CP1}, at which the two conditions hold:

\noindent Condition (i): For any $1 \leq k \leq K_{i}$ and $i \in \mathcal{I}$, we have
\begin{align*}
&	\sum_{k'=1}^k q_{ik'} \ \ = \ \
		1 - \exp \left(- \frac{\omega \delta_{k}^2 \bar{a}_{i}}{\sigma_{i}} \right)
\end{align*}
Condition (ii): For any $1 \leq k \leq K_{i}$ and $(i,j) \in \mathcal{IJ}$, either $\bar{d}_{ijk} = u_{ijk} = 0$, or
\begin{align*}
& \bar{d}_{ijk} \ \ > \ \ 0 \\
& \nu_{k} \bar{d}_{ijk} \ \ln\left(
		\frac{\nu_{k} \bar{d}_{ijk}}{\lambda_{ij} q_{ik} - \nu_{k} \bar{d}_{ijk}} \right)
		\ \ = \ \ \alpha_{ijk} \nu_{k} \bar{d}_{ijk} - \beta_{ij} \nu_{k} u_{ijk}
\end{align*}

\noindent For each $1 \leq k \leq K_{i}$ and $(i,j) \in \mathcal{IJ}$, let
\begin{align*}
x_{ijk} &=
	\begin{cases}
 		u_{ijk} / \bar{d}_{ijk} & \text{ if } \bar{d}_{ijk} \neq 0 \\
 		0 & \text{ otherwise }
	 \end{cases} \\
p_{ijk} &= \frac{ \exp(\alpha_{ijk} - \beta_{ij} x_{ijk}) }
	{ 1 + \exp(\alpha_{ijk} - \beta_{ij} x_{ijk}) }
\end{align*}
Then $(\bar{a}, \bar{d}, \bar{e}, \bar{f}, x, p, q)$ is feasible to \ref{eqn:FP1}, and has the same objective value as $(\bar{a}, \bar{d}, \bar{e}, \bar{f}, u, q)$ in \ref{eqn:CP1}.
\end{lemma}
}

\begin{proof}
Again, the objective value of \ref{eqn:CP1} at $(\bar{a}, \bar{d}, \bar{e}, \bar{f}, u, q)$ is clearly identical to the objective value of \ref{eqn:FP1} at $(\bar{a}, \bar{d}, \bar{e}, \bar{f}, x, p, q)$. Meanwhile, Condition (i) implies that constraint (\ref{eqn:FP1-pickup}) holds in \ref{eqn:FP1}, and Condition (ii) implies that constraints (\ref{eqn:FP1-choice}) and (\ref{eqn:FP1-demand}) hold in \ref{eqn:FP1} (following the algebra in the proof of Lemma \ref{lem:FP1-CP1})). Enough to conclude.
\end{proof}

{\color{purple}
\begin{remark}
We can think of \ref{eqn:CP1} as a relaxation to \ref{eqn:FP1}. Any solution to \ref{eqn:FP1} corresponds to a solution to \ref{eqn:CP1} with the same objective value. However, only some solutions to \ref{eqn:CP1} can be transformed into solutions to \ref{eqn:FP1}. These ``good'' solutions must satisfy the conditions in Lemma \ref{lem:CP1-FP1}. Note the conditions state that the inequality constraints should hold tight that the solution. Thus, the ``good'' solutions are on some ``boundaries'' of  \ref{eqn:CP1}. We will need to prove that any optimal solution must be on those ``boundaries'', and are ``good'' solutions.
\end{remark}
}

{\color{purple}
In short, if the two conditions in Lemma \ref{lem:CP1-FP1} hold, then we can transfer a solution to one problem into a solution to the other with the same objective. If we can show that these two conditions hold at an optimal solution in \ref{eqn:CP1}, then it can be transformed into a feasible solution to \ref{eqn:FP1}, and that solution must be optimal. (Otherwise, we will be able to follow the Lemmas and find a better solution to \ref{eqn:CP1}.) Thus, the theorem can be proved by showing that Condition (i) and Condition (ii) always hold in \ref{eqn:CP1}  at optimality.

Next we show that Condition (ii) in Lemma \ref{lem:CP1-FP1} always hold in \ref{eqn:CP1} at optimality:
}

{\color{blue}
\begin{lemma}
\label{lem:CP1-condition2}
Let $(\bar{a}^*, \bar{d}^*, \bar{e}^*, \bar{f}^*, u^*, q^*)$ be an optimal solution to \ref{eqn:CP1}. Then for any $1 \leq k \leq K_{i}$ and $(i,j) \in \mathcal{IJ}$, either $\bar{d}_{ijk}^* = u_{ijk}^* = 0$, or
\begin{align*}
& \bar{d}_{ijk}^* \ \ > \ \ 0 \\
& \nu_{k} \bar{d}_{ijk}^* \ \ln\left(
		\frac{\nu_{k} \bar{d}_{ijk}^*}{\lambda_{ij} q_{ik}^* - \nu_{k} \bar{d}_{ijk}^*} \right)
		\ \ = \ \ \alpha_{ijk} \nu_{k} \bar{d}_{ijk}^* - \beta_{ij} \nu_{k} u_{ijk}^*
\end{align*}
In other words, Condition (ii) in Lemma \ref{lem:CP1-FP1} holds at $(\bar{a}^*, \bar{d}^*, \bar{e}^*, \bar{f}^*, u^*, q^*)$.
\end{lemma}
}

\begin{proof}
Consider any $1 \leq k \leq K_{i}$ and $(i,j) \in \mathcal{IJ}$. Since $\beta_{ij}>0$, a bigger $u_{ijk}$ is preferred (as we are solving a maximization problem). If $\bar{d}_{ijk}^* = 0$, then we have:
\begin{align*}
& (\lambda_{ij} q_{ik} - \nu_{k} \bar{d}_{ijk}, \
		\nu_{k} \bar{d}_{ijk}, \ \beta_{ij} \nu_{k} u_{ijk} - \alpha_{ijk} \nu_{k} \bar{d}_{ijk}) \\
= \ \ & (\lambda_{ij} q_{ik}, \
		0, \ \beta_{ij} \nu_{k} u_{ijk})
		\ \ \in \ \ \mathcal{K}_{\exp} \\
\implies \quad
& \beta_{ij} \nu_{k} u_{ijk} \ \ \leq \ \ 0
\end{align*}
Note that $u_{ijk}$ is not involved in any other constraint. Thus, if $u_{ijk}^* \neq 0$, then we can find a better feasible solution by setting $u_{ijk} = 0$. That contradicts the fact that $(\bar{a}^*, \bar{d}^*, \bar{e}^*, \bar{f}^*, u^*, q^*)$ is optimal to \ref{eqn:CP1}. Therefore, $u_{ijk}^* = 0$. Similarly, If $\bar{d}_{ijk}^* > 0$, then we have:
\begin{align*}
& (\lambda_{ij} q_{ik} - \nu_{k} \bar{d}_{ijk}, \
		\nu_{k} \bar{d}_{ijk}, \ \beta_{ij} \nu_{k} u_{ijk} - \alpha_{ijk} \nu_{k} \bar{d}_{ijk})
		\ \ \in \ \ \mathcal{K}_{\exp} \\
\implies \quad
& \nu_{k} \bar{d}_{ijk}^* \ \ln\left(
		\frac{\nu_{k} \bar{d}_{ijk}^*}{\lambda_{ij} q_{ik}^* - \nu_{k} \bar{d}_{ijk}^*} \right)
		\ \ \leq \ \ \alpha_{ijk} \nu_{k} \bar{d}_{ijk}^* - \beta_{ij} \nu_{k} u_{ijk}^*
\end{align*}
Again, if the inequality above is not tight, then we can find a better solution to \ref{eqn:CP1} by increasing $u_{ijk}^*$. That contradicts the fact that $(\bar{a}^*, \bar{d}^*, \bar{e}^*, \bar{f}^*, u^*, q^*)$ is optimal to \ref{eqn:CP1}. Therefore, the inequality above must hold at equality. Enough to conclude.
\end{proof}

{\color{purple}
Next we show that Condition (i) in Lemma \ref{lem:CP1-FP1} always hold in \ref{eqn:CP1} at optimality:
}

{\color{blue}
\begin{lemma}
\label{lem:CP1-condition1}
Let $(\bar{a}^*, \bar{d}^*, \bar{e}^*, \bar{f}^*, u^*, q^*)$ be an optimal solution to \ref{eqn:CP1}, such that for any $i \in \mathcal{I}$, if $a_{i}^*>0$, then the following conditions hold:

\noindent Condition (iii): For every $1 \leq k \leq K_{i}$:
\begin{align*}
q_{ik}^* \ \ &> \ \ 0 \\
\sum_{\quad j \in \mathcal{I}} \bar{d}_{ijk}^* \ \ &> \ \ 0	
\end{align*}
\noindent Condition (iv): For every $j \in \mathcal{I}$, either:
\begin{align*}
\frac{\nu_{1} \bar{d}_{ij1}^*}{\lambda_{ij} q_{i1}^*}
		\ \ &> \ \ \cdots
		\ \ > \ \ \frac{\nu_{k} \bar{d}_{ijK_{i}}^*}{\lambda_{ij} q_{iK_{i}}^*}
		\ \ > \ \  0
\end{align*}
or:
\begin{align*}
\bar{d}_{ij1}^*
		\ \ &= \ \ \cdots
		\ \ = \ \ \bar{d}_{ijK_{i}}^*
		\ \ = \ \  0
\end{align*}

\noindent Then, for every $1 \leq k \leq K_{i}$, we have
\begin{align*}
&	\sum_{k'=1}^k q_{ik'}^* \ \ = \ \
		1 - \exp \left(- \frac{\omega \delta_{k}^2 \bar{a}_{i}^*}{\sigma_{i}} \right)
\end{align*}
In other words, Condition (i) in Lemma \ref{lem:CP1-FP1} holds at $(\bar{a}^*, \bar{d}^*, \bar{e}^*, \bar{f}^*, u^*, q^*)$.
\end{lemma}
}

\begin{proof}
Consider any $i \in \mathcal{I}$. If $a_{i}^* = 0$, then since $q^* \geq 0$, we must have $\sum_{k'=1}^k q_{ik'}^* = 0$ for every $1 \leq k \leq K_{i}$. In this case, Condition (i) in Lemma \ref{lem:CP1-FP1} clearly holds for this $i$.

\noindent We assume that $a_{i}^* > 0$ in the rest of this proof. Suppose that
\begin{align*}
&	\sum_{k'=1}^{K_{i}} q_{ik'}^* \ \ < \ \
		1 - \exp \left(- \frac{\omega \delta_{K_{i}}^2 \bar{a}_{i}^*}{\sigma_{i}} \right)
\end{align*}
Then, by increasing $q_{iK_{i}}^*$, we can create some space for some $u_{ijK_{i}}^*$ to be decreased. By decreasing those $u_{ijK_{i}}^*$, we can get a better feasible solution. That contradicts the fact that $(\bar{a}^*, \bar{d}^*, \bar{e}^*, \bar{f}^*, u^*, q^*)$ is optimal. Thus,
\begin{align*}
&	\sum_{k'=1}^{K_{i}} q_{ik'}^* \ \ = \ \
		1 - \exp \left(- \frac{\omega \delta_{K_{i}}^2 \bar{a}_{i}^*}{\sigma_{i}} \right)
\end{align*}

\noindent Next, suppose that
\begin{align*}
\sum_{k=1}^{K_{i}} q_{ik}^* \ \ &= \ \
		1 - \exp \left(- \frac{\omega \delta_{K_{i}}^2 \bar{a}_{i}^*}{\sigma_{i}} \right)  \\
\sum_{k=1}^{K_{i}-1} q_{ik}^* \ \ &< \ \
		1 - \exp \left(- \frac{\omega \delta_{K_{i}-1}^2 \bar{a}_{i}^*}{\sigma_{i}} \right)
\end{align*}
Consider a new solution $(\bar{a}^*, \bar{d}^*, \bar{e}^*, \bar{f}^*, u', q')$, where all variables remain unchanged, except
\begin{align*}
q_{iK_{i}}' \ \ &= \ \
		q_{iK_{i}}^* - \epsilon \\
q_{iK_{i}-1}' \ \ &= \ \
		q_{iK_{i}-1}^* + \epsilon
\end{align*}
and
\begin{align*}
u_{ijK_{i}}' \ \ &= \ \
		\frac{\alpha_{ijK_{i}} \bar{d}_{ijK_{i}}^*}{\beta_{ij}}
		- \frac{\bar{d}_{ijK_{i}}^*}{\beta_{ij}}
 		\ln \left( \frac{\nu_{K_{i}} \bar{d}_{ijK_{i}}^*}
 		{\lambda_{ij} (q_{iK_{i}}^* - \epsilon) - \nu_{K_{i}} d_{ijK_{i}}^*} \right) \\
u_{ijK_{i}-1}' \ \ &= \ \
		\frac{\alpha_{ijK_{i}-1} \bar{d}_{ijK_{i}-1}^*}{\beta_{ij}}
		- \frac{\bar{d}_{ijK_{i}-1}^*}{\beta_{ij}}
		\ln \left( \frac{\nu_{K_{i}-1} \bar{d}_{ijK_{i}-1}^*}
		{\lambda_{ij} (q_{iK_{i}-1}^* + \epsilon) - \nu_{K_{i}-1} d_{ijK_{i}-1}^*} \right)
\end{align*}
for every $j \in \mathcal{I}$ such that
\begin{align*}
\frac{\nu_{1} \bar{d}_{ij1}^*}{\lambda_{ij} q_{i1}^*}
		\ \ &> \ \ \cdots
		\ \ > \ \ \frac{\nu_{k} \bar{d}_{ijK_{i}}^*}{\lambda_{ij} q_{iK_{i}}^*}
		\ \ > \ \  0
\end{align*}
(By Condition (iii), there must be at least one such $j$.)

\noindent We only need to show that, for every such $j \in \mathcal{I}$:
\begin{align*}
\nu_{K_{i}} u_{ijK_{i}}' + \nu_{K_{i}-1} u_{ijK_{i}-1}'
		\ \ > \ \ \nu_{K_{i}} u_{ijK_{i}}^* + \nu_{K_{i}-1} u_{ijK_{i}-1}^*
\end{align*}
with a small enough $\epsilon > 0$. This means $(\bar{a}^*, \bar{d}^*, \bar{e}^*, \bar{f}^*, u', q')$ gives a better objective value, contradicting the fact that  $(\bar{a}^*, \bar{d}^*, \bar{e}^*, \bar{f}^*, u^*, q^*)$ is optimal to \ref{eqn:CP1}.

\noindent Indeed. Let
\begin{align*}
\psi_{j}(\epsilon) \ \ = \ \
\frac{\nu_{K_{i}} \bar{d}_{ijK_{i}}^*}{\beta_{ij}}
 		\ln \left( \frac{\nu_{K_{i}} \bar{d}_{ijK_{i}}^*}
 		{\lambda_{ij} (q_{iK_{i}}^* - \epsilon) - \nu_{K_{i}} \bar{d}_{ijK_{i}}^*} \right)
 + \frac{\nu_{K_{i}-1} d_{ijK_{i}-1}^*}{\beta_{ij}}
		\ln \left( \frac{\nu_{K_{i}-1} \bar{d}_{ijK_{i}-1}^*}
		{\lambda_{ij} (q_{iK_{i}-1}^* + \epsilon) - \nu_{K_{i}-1} \bar{d}_{ijK_{i}-1}^*} \right)
\end{align*}
Note that
\begin{align*}
\nu_{K_{i}} u_{ijK_{i}}' + \nu_{K_{i}-1} u_{ijK_{i}-1}'
		\ \ - \ \ \nu_{K_{i}} u_{ijK_{i}}^* + \nu_{K_{i}-1} u_{ijK_{i}-1}^*
		\ \ = \ \ \psi_{j}(0) - \psi_{j}(\epsilon)
\end{align*}
Thus, it is enough to show that $\psi_{j}'(0) < 0$. We have
\begin{align*}
\psi_{j}'(\epsilon) \ \ = \ \
\frac{\nu_{K_{i}} \bar{d}_{ijK_{i}}^*}{\beta_{ij}}
 		\frac{\lambda_{ij}}{\lambda_{ij} (q_{iK_{i}}^* - \epsilon)
 		- \nu_{K_{i}} \bar{d}_{ijK_{i}}^*}
- \frac{\nu_{K_{i}-1} \bar{d}_{ijK_{i}-1}^*}{\beta_{ij}}
		\frac{\lambda_{ij}}{\lambda_{ij} (q_{iK_{i}-1}^* + \epsilon)
		- \nu_{K_{i}-1} \bar{d}_{ijK_{i}-1}^*}
\end{align*}
That is, $g_{j}'(0) < 0$ holds if and only if
\begin{align*}
\frac{\nu_{K_{i}} \bar{d}_{ijK_{i}}^*}
		{\lambda_{ij} q_{iK_{i}}^* - \nu_{K_{i}} \bar{d}_{ijK_{i}}^*}
		\ \ < \ \
		\frac{\nu_{K_{i}-1} \bar{d}_{ijK_{i}-1}^*}
		{\lambda_{ij} q_{iK_{i}-1}^* - \nu_{K_{i}-1} \bar{d}_{ijK_{i}-1}^*}
\end{align*}
The inequality above holds since
\begin{align*}
\frac{\nu_{K_{i}} \bar{d}_{ijK_{i}}^*}
		{\lambda_{ij} q_{iK_{i}}^*}
		\ \ < \ \
		\frac{\nu_{K_{i}-1} \bar{d}_{ijK_{i}-1}^*}
		{\lambda_{ij} q_{iK_{i}-1}^*}
\end{align*}
and we can conclude that
\begin{align*}
\sum_{k=1}^{K_{i}-1} q_{ik}^* \ \ &= \ \
		1 - \exp \left(- \frac{\omega \delta_{K_{i}-1}^2 \bar{a}_{i}^*}{\sigma_{i}} \right)
\end{align*}
We can repeat this process for $K-2, K-3, \cdots, 1$, and show that
\begin{align*}
&	\sum_{k'=1}^k q_{ik'}^* \ \ = \ \
		1 - \exp \left(- \frac{\omega \delta_{k}^2 \bar{a}_{i}^*}{\sigma_{i}} \right)
\end{align*}
for every $1 \leq k \leq K_{i}$. Enough to conclude.
\end{proof}

{\color{purple}
That is, Condition (i) in Lemma \ref{lem:CP1-FP1} can be ``replaced'' by Condition (iii) and Condition (iv) in Lemma \ref{lem:CP1-condition1}. Next we show that Condition (iii) always hold in \ref{eqn:CP1} at optimality:
}

{\color{blue}
\begin{lemma}
\label{lem:CP1-condition3}
Let $(\bar{a}^*, \bar{d}^*, \bar{e}^*, \bar{f}^*, u^*, q^*)$ be an optimal solution to \ref{eqn:CP1}. If $a_{i}^*>0$, then
\begin{align*}
q_{ik}^* \ \ &> \ \ 0 \\
\sum_{\quad j \in \mathcal{I}} \bar{d}_{ijk}^* \ \ &> \ \ 0	
\end{align*}
for every $1 \leq k \leq K_{i}$. In other words, Condition (iii)  in Lemma \ref{lem:CP1-condition1} holds at $(\bar{a}^*, \bar{d}^*, \bar{e}^*, \bar{f}^*, u^*, q^*)$.
\end{lemma}
}

\begin{proof}
First, given the fact that (1) there is a positive demand, (2) there is a positive number of cars, (3) prices can be arbitrarily large, the optimal solution to \ref{eqn:CP1} must be strictly positive. Thus, there must be some $(i', j') \in \mathcal{IJ}$ and $1 \leq k' \leq K_{i'}$ such that $\bar{d}_{i'j'k'}^*$ is strictly positive.

\noindent Second, suppose that $a_{i}^*>0$, and
\begin{align*}
& \bar{d}_{ij1}^*, \cdots, \bar{d}_{ijK}^* \ \ = \ \ 0	
		\quad , \quad j \in \mathcal{I}
\end{align*}
Then we can decrease $a_{i}^*$ (as well as every $q_{ik}^*$ variables if need), and increase $a_{i'}^*$ by the same amount. By doing so, we are above to slightly increase $q_{i'k'}^*$ and $u_{i'j'k'}^*$ (without changing other variables including $\bar{d}_{i'j'k'}^*$). That will give us a better solution, which contradicts the fact that $(\bar{a}^*, \bar{d}^*, \bar{e}^*, \bar{f}^*, u^*, q^*)$ is optimal.

\noindent Thus, there must be some $j \in \mathcal{I}$ and $1 \leq k \leq K_{i}$ such that $\bar{d}_{ijk}^* > 0$.

\noindent Finally, we show that
\begin{align*}
& \bar{d}_{ijk'}^* \ \ > \ \ 0	
\end{align*}
for every $1 \leq k' \leq K_{i}$. This will complete the proof since it implies
\begin{align*}
q_{ik}^* \ \ &> \ \ 0 \\
\sum_{\quad j \in \mathcal{I}} \bar{d}_{ijk}^* \ \ &> \ \ 0	
\end{align*}

\noindent Indeed. Suppose that $\bar{d}_{ijk'}^* = 0$ for some $1 \leq k' \leq K_{i}$. Consider a new solution $(\bar{a}^*, \bar{d}', \bar{e}^*, \bar{f}^*, u', q')$, where all variables remain unchanged, except
\begin{align*}
q_{ik}' \ \ &= \ \
		q_{ik}^* - \varepsilon \\
q_{ik'}' \ \ &= \ \
		q_{ik'}^* + \varepsilon \\
\bar{d}_{ijk}' \ \ &= \ \ \bar{d}_{ijk}^* - \frac{\epsilon}{\nu_{k}} \\
\bar{d}_{ijk'}' \ \ &= \ \ \frac{\epsilon}{\nu_{k'}}
\end{align*}
and
\begin{align*}
u_{ijk}' \ \ &= \ \
		\frac{\alpha_{ijk} (\bar{d}_{ijk}^* - \epsilon / \nu_{k})}{\beta_{ij}}
		- \frac{\bar{d}_{ijk}^* - \epsilon / \nu_{k}}{\beta_{ij}}
 		\ln \left( \frac{\nu_{k} \bar{d}_{ijk}^* - \epsilon}
 		{\lambda_{ij} (q_{ik}^* - \varepsilon) - (\nu_{k} d_{ijk}^* - \epsilon)} \right) \\
u_{ijk'}' \ \ &= \ \
		\frac{\alpha_{ijk'} \epsilon / \nu_{k'}}{\beta_{ij}}
		- \frac{\epsilon / \nu_{k'}}{\beta_{ij}}
		\ln \left( \frac{\epsilon}
		{\lambda_{ij} (q_{ik'}^* + \varepsilon) - \epsilon} \right)
\end{align*}
for some small enough $\varepsilon > 0$ and $\epsilon > 0$. It is clear that $(\bar{a}^*, \bar{d}', \bar{e}^*, \bar{f}^*, u', q')$ is feasible to \ref{eqn:CP1}. Meanwhile, the increasing rate of $\nu_{k'} u_{ijk'}'$ is infinitely large as $\epsilon$ approaches zero. Thus, we must be able to choose $\varepsilon$ and $\epsilon$ such that $(\bar{a}^*, \bar{d}', \bar{e}^*, \bar{f}^*, u', q')$ gives a better objective value than $(\bar{a}^*, \bar{d}^*, \bar{e}^*, \bar{f}^*, u^*, q^*)$. Enough to conclude.
\end{proof}

{\color{purple}
It left to show that Condition (iv)  in Lemma \ref{lem:CP1-condition1} always hold in \ref{eqn:CP1} at optimality. Indeed:
}

{\color{blue}
\begin{lemma}
\label{lem:CP1-condition4}
Let $(\bar{a}^*, \bar{d}^*, \bar{e}^*, \bar{f}^*, u^*, q^*)$ be an optimal solution to \ref{eqn:CP1}. Consider each $i \in \mathcal{I}$:

\noindent If $a_{i}^*>0$, and for every $1 \leq k \leq K_{i}$:
\begin{align*}
q_{ik}^* \ \ &> \ \ 0
\end{align*}
Then, for every $j \in \mathcal{I}$, either:
\begin{align*}
\frac{\nu_{1} \bar{d}_{ij1}^*}{\lambda_{ij} q_{i1}^*}
		\ \ &> \ \ \cdots
		\ \ > \ \ \frac{\nu_{k} \bar{d}_{ijK_{i}}^*}{\lambda_{ij} q_{iK_{i}}^*}
		\ \ > \ \  0
\end{align*}
or:
\begin{align*}
\bar{d}_{ij1}^*
		\ \ &= \ \ \cdots
		\ \ = \ \ \bar{d}_{ijK_{i}}^*
		\ \ = \ \  0
\end{align*}
In other words, Condition (iii) in Lemma \ref{lem:CP1-condition1} implies Condition (iv).
\end{lemma}
}

\begin{proof}
Consider any $(i,j) \in \mathcal{IJ}$. With an identical approach in the proof of \ref{lem:CP1-condition3}, we can show that, if $\bar{d}_{ijk} > 0$ for any $1 \leq k \leq K_{i}$, then $\bar{d}_{ijk} > 0$ for every $1 \leq k \leq K_{i}$. In the rest of this proof, we assume that $\bar{d}_{ijk} > 0$ for every $1 \leq k \leq K_{i}$.

\noindent We prove this Lemma by contradiction. Specifically, we show that if for any $2 \leq k \leq K_{i}$:
\begin{align*}
\frac{\nu_{k-1} \bar{d}_{ijk-1}^*}{\lambda_{ij} q_{ik-1}^*}
\ \ \leq \ \
\frac{\nu_{k} \bar{d}_{ijk}^*}{\lambda_{ij} q_{ik}^*}
\end{align*}
Then the current solution $(\bar{a}^*, \bar{d}^*, \bar{e}^*, \bar{f}^*, u^*, q^*)$ is not optimal, and we can find a better solution. Consider solution $(\bar{a}', \bar{d}', \bar{e}^*, \bar{f}^*, u^*, q')$, where all decision variables remain unchanged, except
\begin{align*}
\bar{a}_{i}' \ \ &= \ \
		\bar{a}_{i}^* + \frac{\epsilon}{\nu_{k}} - \frac{\epsilon}{\nu_{k-1}}
		\ \ \geq \ \ \bar{a}_{i}^*
\end{align*}
and
\begin{align*}
\bar{d}_{ijk}' \ \ &= \ \ \bar{d}_{ijk}^* - \frac{\epsilon}{\nu_{k}} \\
\bar{d}_{ijk-1}' \ \ &= \ \ \bar{d}_{ijk-1}^* + \frac{\epsilon}{\nu_{k-1}} \\
u_{ijk}' \ \ &= \ \ \frac{\alpha_{ijk} (\bar{d}_{ijk}^* - \epsilon / \nu_{k})}{\beta_{ij}}
		- (\frac{\bar{d}_{ijk}^* - \epsilon / \nu_{k}}{\beta_{ij}})
 		\ln( \frac{\nu_{k} \bar{d}_{ijk}^* - \epsilon}
 				{\lambda_{ij} q_{ik}^* - (\nu_{k} \bar{d}_{ijk}^* - \epsilon)} ) \\
u_{ijk-1}' \ \ &= \ \ \frac{\alpha_{ijk-1} (\bar{d}_{ijk-1}^* + \epsilon / \nu_{k-1})}{\beta_{ij}}
		-  (\frac{\bar{d}_{ijk-1}^* + \epsilon / \nu_{k-1}}{\beta_{ij}})
		\ln( \frac{\nu_{k-1} \bar{d}_{ijk-1}^* + \epsilon}
				{\lambda_{ij} q_{ik-1}^* - (\nu_{k-1} \bar{d}_{ijk-1}^* + \epsilon)} )
\end{align*}
for every $j \in \mathcal{I}$ and some small enough $\epsilon > 0$. Note that
\begin{align*}
\nu_{k} \bar{d}_{ijk}' + \nu_{k-1} \bar{d}_{ijk-1}'
		\ \ = \ \ \nu_{k} \bar{d}_{ijk}^* + \nu_{k-1} \bar{d}_{ijk-1}^*
		\quad, \quad j \in \mathcal{I}
\end{align*}
Thus, all other flows in the network are not affected, which means all other constraints still hold. Meanwhile, since $\phi_{ijk} \geq \phi_{ijk-1}$, we have
\begin{align*}
\nu_{k} \phi_{ijk} \bar{d}_{ijk}' + \nu_{k-1} \phi_{ijk-1} \bar{d}_{ijk-1}'
		\ \ \leq \ \ \nu_{k} \phi_{ijk} \bar{d}_{ijk}^* + \nu_{k-1} \phi_{ijk-1} \bar{d}_{ijk-1}^*
		\quad, \quad j \in \mathcal{I}
\end{align*}
Thus, we only need to show that
\begin{align*}
\nu_{k} u_{ijk}' + \nu_{k-1} u_{ijk-1}' \ \ > \ \ \nu_{k} u_{ijk}^* + \nu_{k-1} u_{ijk-1}^*
		\quad, \quad j \in \mathcal{I}
\end{align*}
Indeed. Let
\begin{align*}
\varrho_j(\epsilon) \ \ \defi \ \
		& \nu_{k} u_{ijk}' + \nu_{k-1} u_{ijk-1}' \\
= \ \ &\frac{\alpha_{ijk} (\nu_{k} \bar{d}_{ijk}^* - \epsilon)}{\beta_{ij}}
		- (\frac{\nu_{k} \bar{d}_{ijk}^* - \epsilon}{\beta_{ij}})
 				\ln( \frac{\nu_{k} \bar{d}_{ijk}^* - \epsilon}{\lambda_{ij} q_{ik}^* - (\nu_{k} \bar{d}_{ijk}^* - \epsilon)} ) \\
+ &\frac{\alpha_{ijk-1} (\nu_{k-1} \bar{d}_{ijk-1}^* + \epsilon)}{\beta_{ij}}
		-  (\frac{\nu_{k-1} \bar{d}_{ijk-1}^* + \epsilon}{\beta_{ij}})
		\ln( \frac{\nu_{k-1} \bar{d}_{ijk-1}^* + \epsilon}{\lambda_{ij} q_{ik-1}^* - (\nu_{k-1} \bar{d}_{ijk-1}^* + \epsilon)} )
\end{align*}
We have
\begin{align*}
\varrho_j'(\epsilon) \ \ = \ \
&- \frac{\alpha_{ijk}}{\beta_{ij}} + \frac{1}{\beta_{ij}} \Bigg(
 		\ln( \frac{\nu_{k} \bar{d}_{ijk}^*}{\lambda_{ij} q_{ik}^* - \nu_{k} \bar{d}_{ijk}^*} )
		+ \frac{\lambda_{ij} q_{ik}^*}{\lambda_{ij} q_{ik}^* - \nu_{k} \bar{d}_{ijk}^*} \Bigg) \\
&+ \frac{\alpha_{ijk-1}}{\beta_{ij}} - \frac{1}{\beta_{ij}} \Bigg(
 		\ln( \frac{\nu_{k-1} \bar{d}_{ijk-1}^*}{\lambda_{ij} q_{ik-1}^* - \nu_{k-1} \bar{d}_{ijk-1}^*} )
		+ \frac{\lambda_{ij} q_{ik-1}^*}{\lambda_{ij} q_{ik-1}^* - \nu_{k-1} \bar{d}_{ijk-1}^*} \Bigg)
\end{align*}
Recall that $\ln(x/(1-x)) + 1/(1-x)$ is an increasing function on $(0,1)$. Since $\alpha_{ijk-1} > \alpha_{ijk}$ and
\begin{align*}
\frac{\nu_{k-1} \bar{d}_{ijk-1}^*}{\lambda_{ij} q_{ik-1}^*}
\ \ \leq \ \
\frac{\nu_{k} \bar{d}_{ijk}^*}{\lambda_{ij} q_{ik}^*}
\end{align*}
we know that $\varrho_j'(\epsilon) > 0$. Note that
\begin{align*}
\varrho_j(0) \ \ \defi \ \
		& \nu_{k} u_{ijk}^* + \nu_{k-1} u_{ijk-1}^*
\end{align*}
Thus, $(\bar{a}', \bar{d}', \bar{e}^*, \bar{f}^*, u^*, q')$ is a better solution. Enough to conclude. \\
\end{proof}

\newpage


\subsection{Proof of Theorem \ref{thm:CP2-FP2}}

{\color{purple}
Theorem \ref{thm:CP2-FP2} can be proved in a very similar way as Theorem \ref{thm:CP1-FP1} is proved:
}

First, let $(\bar{a}, \bar{e}, \bar{f}, x, p, q)$ be a feasible solution to \ref{eqn:FP2}. Let $u_{ij} = \bar{f}_{ij} x_{ij}$ for each $(i,j) \in \mathcal{IJ}$. Then $(\bar{a}, \bar{e}, \bar{f}, u, q)$ is a feasible solution to \ref{eqn:CP2} with identical objective value (same as Lemma \ref{lem:FP1-CP1}).

Second, let $(\bar{a}^*, \bar{e}^*, \bar{f}^*, u^*, q^*)$ be an optimal solution to \ref{eqn:CP2}. Consider each $(i,j) \in \mathcal{IJ}$. We know that: If $\bar{f}_{ij}^* = 0$, then $u_{ij}^* = 0$. If $\bar{f}_{ij}^* > 0$, then
\begin{align*}
& \mu_{ij} \bar{f}_{ij}^* \ \ln\left(
		\frac{\mu_{ij} \bar{f}_{ij}^*}{\lambda_{ij} - \mu_{ij} \bar{f}_{ij}^*} \right)
		\ \ = \ \ \alpha_{ij0} \mu_{ij} \bar{f}_{ij}^* - \beta_{ij} \mu_{ij} u_{ij}^*
\end{align*}
(Otherwise, we can increase $u_{ij}^*$ and get a better solution, same as in Lemma  \ref{lem:CP1-condition2}.)

Third, the optimal objective value must be strictly positive. To show this, we can construct a feasible solution with positive objective value as follow: (1) Recall that $\mathcal{IJ} \neq \varnothing$. We can set $q_i = 1$ for some $i$ such that $(i,j) \in \mathcal{IJ}$ for some $j$; (2) Recall that $\lambda_{ij} > 0$ for every $j$ such that $(i,j) \in \mathcal{IJ}$. We can set $\bar{f}_{ij}$ equal to a $\epsilon > 0$ for every $j$ such that $(i,j) \in \mathcal{IJ}$ (all other $\bar{f}$ variables are set to $0$); (3) Recall that between any two zones, there must be a sequence of paths that empty cars can flow through. Thus we can increase some of the $\bar{e}$ variables, such that the network flow is balanced (all other $\bar{e}$ variables are set to $0$); (4) As $\bar{f}_{ij}$ ($\epsilon$) approaches zero, the gradient of
\begin{align*}
\mu_{ij} \bar{f}_{ij}^* \ \ln\left(
		\frac{\mu_{ij} \bar{f}_{ij}^*}{\lambda_{ij} - \mu_{ij} \bar{f}_{ij}^*} \right)
\end{align*}
goes to infinite. Thus, there must be some small enough $\epsilon$ such that the network flow is balanced, the cardinality constraint is satisfied (by controlling the $\bar{a}$ variables), and the objective value is strictly positive (since $u_{ij}$ can increase "infinitely fast" as a function of $\epsilon$).

The result above implies $\bar{f}^* > 0$. Indeed. Suppose there is some $i$ where $\bar{f}_{ij} = 0$ for some $j$ such that $(i,j) \in \mathcal{IJ}$. Then we can construct a better solution by: (1) Set $q_i = 1$, (2) Scale down all car flows in the network by a factor of $1 - \epsilon$, where $\epsilon > 0$ can be arbitrarily close to $0$; (3) For every $j$ such that $(i,j) \in \mathcal{IJ}$, increase $\bar{f}_{ij}$ by a similar amount, and then increase the corresponding empty car flows to keep the network balanced, until the cardinality constraint is satisfied. As $\epsilon$ approaches zero, the increase rate in profit goes to infinity.

Note that $\bar{f}^* > 0$ then implies $q^* = 1$. Otherwise, we will be able to get a better solution by increasing some $q_i^*$, and increasing the corresponding $u_{ij}^*$ variables without changing anything else.

Now, For each $(i,j) \in \mathcal{IJ}$, let
\begin{align*}
x_{ij}^* \ \ &= \ \
		\begin{cases}
 			u_{ij}^* / \bar{f}_{ij}^* & \text{ if } \bar{f}_{ij}^* > 0 \\
 			0 & \text{ otherwise }
	 	\end{cases} \\
p_{ij}^* \ \ &= \ \
		\frac{ \exp(\alpha_{ij0} - \beta_{ij} x_{ij}^*) }
		{ 1 + \exp(\alpha_{ij0} - \beta_{ij} x_{ij}^*) }
\end{align*}
From the arguments above, we know that $(\bar{a}^*, \bar{e}^*, \bar{f}^*, x^*, p^*, q^*)$ is a feasible solution to \ref{eqn:FP2} with identical objective value (same as Lemma \ref{lem:CP1-FP1}). Combining with the first argument, we know that $(\bar{a}^*, \bar{e}^*, \bar{f}^*, x^*, p^*, q^*)$ is an optimal solution to \ref{eqn:FP2}.

\newpage


\subsection{Proof of Theorem \ref{thm:CP3-FP3}}

{\color{purple}
The proof of Theorem \ref{thm:CP3-FP3} is identical to the proof of Theorem \ref{thm:CP1-FP1}, except we need to replace each $\bar{a}_{i}$ by $\bar{a}_{i} + \sum_{j \in \widetilde{\mathcal{J}}_{i}} \zeta_{ji} \bar{e}_{ji}$.
}


\end{document}